\PassOptionsToPackage{numbers,sort&compress}{natbib}
\documentclass{article}

\usepackage[preprint]{neurips_2024}
\usepackage[utf8]{inputenc} %
\usepackage[T1]{fontenc}    %
\usepackage{hyperref}       %
\usepackage{url}            %
\usepackage{booktabs}       %
\usepackage{amsfonts}       %
\usepackage{nicefrac}       %
\usepackage{microtype}      %
\usepackage{xcolor}         %

\usepackage{amsmath,algorithm}
\usepackage{algpseudocode}
\usepackage{stfloats}
\usepackage{graphicx}
\usepackage{subfigure}
\usepackage{booktabs}
\usepackage{amsmath}
\usepackage{amssymb}
\usepackage{mathtools}
\usepackage{amsthm}
\usepackage[capitalize,noabbrev]{cleveref}
\usepackage[textsize=tiny]{todonotes}

\theoremstyle{plain}
\newtheorem{theorem}{Theorem}[section]
\newtheorem{proposition}[theorem]{Proposition}
\newtheorem{lemma}[theorem]{Lemma}

\theoremstyle{definition}

\newtheorem{example}[theorem]{Example}
\theoremstyle{remark}
\newtheorem{remark}[theorem]{Remark}

\theoremstyle{plain}
\newtheorem{innernamedtheorem}{Theorem}[section]
\newenvironment{named}[1]{%
    \def\theoremname{#1}%
    \begin{innernamedtheorem}[\theoremname]%
}{%
    \end{innernamedtheorem}%
}

\title{
    S-SOS: Stochastic Sum-Of-Squares
    for Parametric Polynomial Optimization
}

\author{
  Richard L. Zhu \\
  Department of Computational and Applied Mathematics\\
  University of Chicago\\
  Chicago, IL 60637 \\
  \texttt{richardzhu@uchicago.edu} \\
  \AND
  Mathias Oster \\
  Institute of Geometry and Practical Mathematics \\
  RWTH Aachen \\
  Aachen, Germany \\
  \texttt{oster@igpm.rwth-aachen.de} \\
  \AND
  Yuehaw Khoo \\
  Department of Statistics \\
  University of Chicago\\
  Chicago, IL 60637 \\
  \texttt{yuehaw.khoo@uchicago.edu}
}

\begin{document}

\maketitle

\begin{abstract}
Global polynomial optimization is an important tool across applied mathematics, with many applications in operations research, engineering, and physical sciences.
In various settings, the polynomials depend on external parameters that may be random.
We discuss a stochastic sum-of-squares (S-SOS) algorithm based on the sum-of-squares hierarchy that constructs a series of semidefinite programs to jointly find strict
lower bounds on the global minimum and extract candidates for parameterized global minimizers.
We prove quantitative convergence of the hierarchy as the degree increases and use it to solve unconstrained and constrained polynomial optimization problems parameterized by random variables.
By employing $n$-body priors from condensed matter physics to induce sparsity, we can use S-SOS to produce solutions and uncertainty intervals for sensor network localization problems
containing up to 40 variables and semidefinite matrix sizes surpassing $800 \times 800$.
\end{abstract}

\section{Introduction}\label{main:intro}

Many effective nonlinear and nonconvex optimization techniques use local information to identify local minima. But it is often the case that we want to find global optima. Sum-of-squares (SOS) optimization is a powerful and general technique in this setting.

The core idea is as follows: suppose we are given polynomials $g_1, \ldots, g_m, f$ where each function is on
$\mathbb{R}^n \to \mathbb{R}$ and we seek to determine the minimum value of $f$ on the closed set $\mathcal{S}$: $\mathcal{S} = \{ x \in \mathbb{R}^n \; | \; g_i(x) \geq 0 \; \forall \; i = 1, \ldots, m \}$.
Our optimization problem is then to find $ \inf_{x \in \mathbb{R}^n} \{ f(x) | x \in \mathcal{S} \} $.
An equivalent formulation is to find the largest constant $c \in \mathbb{R}$ (i.e. the tightest lower bound)
that can be subtracted from $f$ such that $f - c \geq 0$ over the set $\mathcal{S}$.
This reduction converts a polynomial optimization problem over a semialgebraic set to the problem of checking polynomial non-negativity. This problem is NP-hard in general \cite{garey_computers_2009}, therefore one instead resorts to checking if $f - c$ is a sum-of-squares (SOS) function, e.g. in the unconstrained setting
where $\mathcal{S} = \mathbb{R}^n$ one seeks to find some polynomials $h_k: \mathbb{R}^n \to \mathbb{R}$ such that $ f - c = \sum_k h_k^2 $. If such a decomposition can be found, then we have an easily checkable certification
that $f - c \geq 0$, as all sum-of-squares are non-negative but not all non-negative
functions are sum-of-squares.

Notably, if we restrict the $h_k$ to have maximum degree $d$, the search for a degree-$2d$ SOS decomposition of a function can be automated as a semidefinite program (SDP) \cite{nesterov_squared_2000,lasserre_global_2001,laurent_sums_2009}. Solving this SDP for varying degrees $d$ generates the well-known
Lasserre (SOS) hierarchy. A given degree $d$ corresponds to
a particular level of the hierarchy.
Solving this SDP produces a lower bound $c_d$ which has been proven to converge to the true global minimum $c^* = \inf_x f(x)$ as $d$ increases,
with finite convergence
($c_d = c^*$ at finite $d$) for functions with second-order local optimality conditions \cite{nie_optimality_2014,bach_exponential_2023} and asymptotic
convergence with milder assumptions thanks to representation theorems for positive polynomials from real algebraic geometry \cite{putinar_positive_1993,schmudgen_moment_2017}.
Further work has elucidated both theoretical 
implications \cite{putinar_positive_1993, lasserre_global_2001, lasserre_moment-sos_2018, lasserre_moment-sos_2023} and useful applications 
of SOS to disparate fields
\cite{parrilo_structured_2000, nie_sum_2009, de_klerk_complexity_2008, nie_optimality_2014,bach_exponential_2023,ahmadi_dsos_2019, papp_sum--squares_2019} (see further discussion in \cref{appendix:related-work}).

Motivated by the sum-of-squares certification for a lower bound $c$ on a function $f(x)$,
we generalize to the case where the function to be minimized has additional parameters, i.e. $f(x, \omega)$ where $x$ are variables
and $\omega$ are parameters drawn
from some probability distribution $\omega \sim \nu(\omega)$. We seek a function $c(\omega)$ that is the tightest lower bound to $f(x, \omega)$ everywhere: $f(x, \omega) \geq c(\omega)$ with $c(\omega) \to \inf_x f(x, \omega)$. This setting was originally presented in \cite{lasserre_jointmarginal_2010} as a ``Joint and Marginal'' approach to parametric polynomial optimization.
With the view that $\omega \sim \nu(\omega)$ and seeking to parameterize the minimizers $x^*(\omega) = \text{argmin}_x f(x, \omega)$,
we are reminded of some of the prior work in polynomial chaos, where a system of stochastic variables is expanded into a deterministic function of those stochastic variables \cite{sudret_global_2008,najm_uncertainty_2009}.

\textbf{Contributions and outline.}
Our primary contributions are a quantitative convergence proof for the Stochastic Sum-of-Squares (S-SOS) hierarchy of semidefinite programs (SDPs), a formulation of a new hierarchy (the cluster basis hierarchy) that uses the structure of a problem to sparsify the SDP, and numerical results on its application to the sensor network localization problem.

In \cref{section:ssos}, we review the S-SOS hierarchy of SDPs \cite{lasserre_jointmarginal_2010} and its primal and dual formulations (\cref{section:ssos:formulation}). We then detail how different hierarchies can be constructed (\cref{section:ssos:variations:lasserre-hierarchy}).
    Finally, in \cref{subsection:convergence}
    (complete proofs in \cref{appendix:convergence-proof}) we specialize to compact $X \times \Omega$ and outline the proofs for two theorems on quantitative convergence
    (the gap between the optimal values of the degree-$2s$ S-SOS SDP and the ``tightest lower-bounding'' optimization problem goes $\to 0$ as $s \to \infty$)
    of the S-SOS hierarchy for trigonometric polynomials on $[0, 1]^{n} \times [0, 1]^d$ following the kernel formalism of
    \cite{fang_sum--squares_2021,bach_exponential_2023,slot_sum--squares_2023}.
    The first one applies in the general case and the second one applies to the case where $d=1$.

In \cref{section:numerics} we review the hierarchy's 
    applications in parametric polynomial minimization and uncertainty quantification,
    focusing on several variants of sensor network localization on $X \times \Omega = [-1, 1]^n \times [-1, 1]^d$. We present numerical results for the accuracy of the extracted solutions that result from S-SOS, comparing to other approaches to parametric polynomial optimization,
    including a simple Monte Carlo-based method.

\section{Stochastic Sum-of-squares (S-SOS)}\label{section:ssos}

\subsubsection{Notation}\label{section:ssos:notation}

Let $\mathcal{P}(S)$ be the space of polynomials on $S$, where $S \in \{X, \Omega\}$.
$X \subseteq \mathbb{R}^n$ and $\Omega \subseteq \mathbb{R}^d$, respectively,
where $X$ and $\Omega$ are (not-necessarily compact) subsets of their respective ambient spaces $\mathbb{R}^n$ and $\mathbb{R}^d$.
A polynomial in $\mathcal{P}(X)$ can be written as $p(x) = \sum_{\alpha \in \mathbb{Z}^n_{\geq 0}} c_{\alpha} x^{\alpha} \in \mathcal{P}(X)$
(substituting $n \to d, x \to \omega, X \to \Omega$ for a polynomial in $\Omega$).
Let $x := (x_1, \ldots, x_n), \omega := (\omega_1, \ldots, \omega_d)$, $\alpha$ be a multi-index (size given by context), and $c_{\alpha}$
be the polynomial coefficients.
Let $\mathcal{P}^s(S)$ for some $s \in \mathbb{Z}_{\geq 0}, S \in \{X, \Omega\}$ denote the subspace of $\mathcal{P}(S)$
consisting of polynomials of degree $\leq s$, i.e. polynomials where the multi-indices of the monomial terms satisfy $||\alpha||_{1} \leq s$.
$\mathcal{P}_{\text{SOS}}(X \times \Omega)$ refers to the space of polynomials on $X \times \Omega$ that can be expressible as a sum-of-squares in $x$ and $\omega$ jointly,
and $\mathcal{P}_{\text{SOS}}^s(X \times \Omega)$ be the same space restricted to polynomials of degree $\leq s$.
Additionally, $W \succcurlyeq 0$ for a matrix $W$ denotes that $W$ is symmetric positive semidefinite (PSD). Finally, $\mathbb P(\Omega)$ denotes the set of Lebesgue probability measures on $\Omega$.
For more details, see \cref{appendix:notation}.

\subsection{Formulation of S-SOS hierarchy}\label{section:ssos:formulation}

We present two formulations of the 
S-SOS hierarchy that are dual to each other in the sense of Fenchel duality \cite{Rockafellar,boyd_convex_2004}.
The primal problem seeks to find
the tightest lower-bounding function and the dual problem
seeks to find a minimizing probability distribution.
Note that the ``tightest lower bound'' approach is dual to the 
``minimizing distribution'' approach, otherwise known as a ``joint and marginal'' moment-based approach
originally detailed in \cite{lasserre_jointmarginal_2010}.

\subsubsection{Primal S-SOS: The tightest lower-bounding function}

Consider a polynomial $f(x, \omega): \mathbb{R}^{n + d} \to \mathbb{R}$ with
$x \in X \subseteq \mathbb{R}^n, \omega \in \Omega \subseteq \mathbb{R}^d$
equipped with a probability measure $\nu(\omega)$.
We interpret $x$ as our optimization variables and $\omega$ as noise
parameters, and seek a lower-bounding function
$c^*(\omega)$ such that $f(x, \omega) \geq c^*(\omega)$ for all $x, \omega$.
In particular, we want the tightest lower bound $c^*(\omega) = \inf_{x \in X} f(x, \omega)$.
Note that even when $f(x, \omega)$ is polynomial,
the tightest lower bound $c^*(\omega)$ can be non-polynomial.
A simple example is the function $f(x, \omega) = (x - \omega)^2 + (\omega x)^2$, which has $c^*(\omega) = \inf_x f(x, \omega) = \omega^4 / (1 + \omega^2)$ (\cref{appendix:simple-potential:analytic-solution}).

For us to select the ``best'' lower-bounding function, we want to maximize the expectation
of the lower-bounding function $c(\omega)$
under $\omega \sim \nu(\omega)$ while
requiring $f(x, \omega) - c(\omega) \geq 0$,
giving us the following optimization problem
over $L^1$-integrable lower-bounding functions:
\begin{align}
p^* = \sup_{c \in L^1(\Omega)}        \quad & \int c(\omega) \mathrm{d}\nu(\omega) \label{eq:primal-sos} \\
                   \text{s.t.} \quad & f(x, \omega) - c(\omega) \geq 0 \notag
\end{align}

Even if we restricted $c(\omega)$ to be 
polynomial so that the residual $f(x, \omega) - c(\omega)$ is also polynomial, we would still have a 
challenging nonconvex optimization problem over non-negative polynomials. In SOS optimization, we take a relaxation and require the residual to be SOS:
$f(x, \omega) - c(\omega) \in \mathcal{P}_{\text{SOS}}(X \times \Omega)$.
Doing the SOS relaxation of the non-negative \cref{eq:primal-sos} and restricting $c(\omega)$, i.e. $f(x,\omega) - c(\omega)$ to polynomials
of degree $\leq 2s$ gives us \cref{eq:primal-sos-deg}, which we call the primal S-SOS degree-$2s$ SDP:
\begin{align}
p^*_{2s} = \sup_{c\in\mathcal{P}^{2s}(\Omega),W \succcurlyeq 0}        \quad & \int c(\omega) \mathrm{d}\nu(\omega) \label{eq:primal-sos-deg}     \\
    \text{s.t.} \quad &  f(x, \omega) - c(\omega) = m_s(x, \omega)^T W m_s(x, \omega) \notag
\end{align}
where $m_s(x, \omega)$ is a basis function $X \times \Omega \to \mathbb{R}^{a(n,d,s)}$
containing monomial terms of degree $\leq s$ written as a column vector,
and $W \in \mathbb{R}^{a(n,d,s) \times a(n,d,s)}$ a symmetric PSD matrix. Here, $a(n,d,s)$
represents the dimension of the basis function, which depends on the degree $s$
and on the dimensions $n, d$.
For this formulation to find the best degree-$2s$ approximation to the lower-bounding function,
we require
$ g(x, \omega) = m_s(x, \omega)^T W m_s(x, \omega) $
to span $\mathcal{P}^{2s}(X \times \Omega)$.
Selecting all combinations of standard monomial terms of degree $\leq s$ suffices and 
results in a basis function with size $a(n,d,s) = \binom{n+d+s}{s}$.

\subsubsection{Dual S-SOS: A minimizing distribution}

The formal dual to \cref{eq:primal-sos} (proof of duality in \cref{appendix:formal-dual}) seeks to find a ``minimizing distribution'' $\mu(x, \omega)$, i.e. a probability distribution
that places weight on the minimizers of $f(x, \omega)$
subject to the constraint that the marginal $\mu_X(\omega)$ matches $\nu(\omega)$:
\begin{align}
d^* = \inf_{\mu \in \mathbb P(X \times\Omega)}     \quad & \int f(x, \omega) \mathrm{d}\mu(x, \omega) \label{eq:dual-sdp} \\
    \text{s.t.} \quad & \int_X \mathrm{d}\mu(x, \omega) = \mu_X(\omega) = \nu(\omega) \notag
\end{align}
where we have written $\mathbb P(X \times \Omega)$ as the space of joint probability distributions on $X \times \Omega$
and $\mu_X(\omega)$ is the marginal of $\mu(x, \omega)$ with respect to $\omega$, obtained via disintegration.

For the primal, we considered polynomials of degree $\leq 2s$.
We do the same here.
The formal dual becomes a tractable SDP, where the objective turns into
moment-minimization and the constraints become moment-matching.
Following \cite{lasserre_global_2001,nie_sum_2009}, let $M \in \mathbb{R}^{a(n,d,s) \times a(n,d,s)}$
be the symmetric PSD moment matrix
with entries defined as
$ M_{i, j} = \int_{X \times \Omega} m_s^{(i)}(x, \omega) m_s^{(j)}(x, \omega) d\mu(x, \omega) $
where $m_s^{(i)}(x, \omega)$ is the $i$-th element of the basis function $m_s$.
Let $y \in \mathbb{R}^{b(n,d,s)}$ be the moment vector of independent moments that completely specifies $M$,
e.g. in the case that we use all standard monomials of degree $\leq s$ and have $a(n, d, s) = \binom{n+d+s}{s}$, then $b(n,d,s) = \binom{n+d+2s}{2s}$.
We write $M(y)$ as
the moment matrix that is formed from these independent moments. We have
$ y_{\alpha(i, j)} = \int_{X \times \Omega} m_{s}^{(i)}(x, \omega) m_s^{(j)}(x, \omega) d\mu(x, \omega) $
where the multi-index $\alpha(i, j) \in \mathbb{Z}_{\geq 0}^{n+d}$ corresponds to
the sum of the multi-indices corresponding to the 
$i$-th entry and the $j$-th entry of $m_s(x, \omega)$.

We write $f(x, \omega)$ in terms of the monomials $f(x, \omega) = \sum_{||\alpha||_{1} \leq 2s} f_{\alpha} [x, \omega]^{\alpha}$,
where $[x, \omega]$ is the concatenation of the $n+d$ variables from $x, \omega$ and
$\alpha \in \mathbb{Z}_{\geq 0}^{n+d}$ is a multi-index.
Note that every monomial $[x, \omega]^{\alpha}$ has a corresponding moment $y_{\alpha}$: $ \int [x, \omega]^{\alpha} \mathrm{d}\mu(x,\omega) = y_{\alpha}$.
We then observe that the integral in the objective reduces to a dot product between
the coefficients of $f$ and the moment vector:
$$ \int f(x, \omega) \mathrm{d}\mu(x, \omega) = \int \sum_{\alpha} f_{\alpha} [x, \omega]^{\alpha} \mathrm{d}\mu(x, \omega) = \sum_{\alpha} f_{\alpha} y_{\alpha} $$
After converting the distribution-matching constraint $\mu_X(\omega) = \nu(\omega)$ in \eqref{eq:dual-sdp} into equality constraints
on the moments of $\omega$ up to degree $2s$,
we obtain the following dual S-SOS degree-$2s$ SDP:
\begin{align}
d^*_{2s} = \inf_{y \in \mathbb{R}^{b(n,d,s)}} \quad & \sum_{||\alpha||_{1} \leq 2s} f_{\alpha} y_{\alpha} \label{eq:dual-sdp-deg} \\
    \text{s.t.} \quad & M(y) \succcurlyeq 0 \notag \\
    & y_{\alpha} = m_{\alpha} \; \forall \; (\alpha, m_{\alpha}) \in \mathcal{M}_{\nu} \notag
\end{align}

We write $\mathcal{M}_{\nu}$ as the set of $(\alpha, m_{\alpha})$ representing the moment-matching constraints on $\omega^{\alpha}$ up to degree-$2s$,
i.e. we want to set $\int_{X \times \Omega} \omega^{\alpha} \mathrm{d}\mu(x, \omega) = \int_{\Omega} \omega^{\alpha} \mathrm{d}\nu(\omega) = m_{\alpha}$
for all multi-indices $\alpha \in \mathbb{Z}_{\geq 0}^{d}$ with $||\alpha||_1 \leq 2s$.
There are $\binom{d + 2s}{2s}$ multi-indices $\alpha \in \mathbb{Z}^{n+d}_{\geq 0}, ||\alpha||_1 \leq 2s$ where only the $d$ entries associated with $\omega$ are non-zero, and therefore the number of moment-matching constraints is $|\mathcal{M}_{\nu}| = \binom{d + 2s}{2s}$.
Note that the moment matrix $M(y) \in \mathbb{R}^{a(n,d,s) \times a(n,d,s)}$ is a symmetric PSD matrix and is the dual variable to the primal $W$.
Observe also that we require the moments of $\nu(\omega)$ of degree up to $2s$ to be bounded.
\eqref{eq:dual-sdp-deg} is often a more convenient
form than \eqref{eq:primal-sos-deg}, especially when working with additional equality
or inequality constraints, as we will see in \cref{section:numerics}.
For concrete examples of the primal and dual SDPs with explicit constraints,
see \cref{appendix:sdp-examples}.

\subsection{Variations}

In this section, we detail two ways of building a hierarchy, one based
on the maximum degree of monomial terms in the basis
function (Lasserre) and a novel one based on the 
maximum number of interactions occurring in the terms of the basis function (cluster basis).
To define any SOS hierarchy, we first select a monomial basis. Some examples include the standard monomial basis $x_1, \ldots, x_n$, trigonometric/Fourier 1-periodic monomial basis $\sin x_1, \cos x_1, \ldots, \sin x_n, \cos x_n$), or others. Using this basis, we write down a basis function $m(x)$ which comprises some combinations of monomials.
Squared linear combinations of the basis functions then span a
SOS space of functions: $\mathcal{H}: \{ (\sum_i h_i m_i(x))^2 \}$.

\subsubsection{Standard Lasserre hierarchy}\label{section:ssos:variations:lasserre-hierarchy}

In the Lasserre hierarchy, the basis function $m_s(x)$ 
is composed of all combinations of monomials 
up to degree $s \in \mathbb{Z}_{> 0}$ and a given level of the hierarchy is set by the maximum degree $s$.
The basis function consists of terms $x^{\alpha}$ with $\alpha$ a multi-index and $||\alpha||_1 \leq s$.
The degree-$2s$ SOS function space parameterized by this basis function is that spanned by $m_s(x)^T W m_s(x)$ for PSD $W$, i.e. the functions that 
can result from squaring any linear combination of degree-$s$ polynomials that can be generated 
from our basis $m_s(x)$.
As we increase the degree $s$, our basis function gets
larger and our S-SOS SDP objective values converge to the optimal value of
the ``tightest lower-bounding'' problem \cref{eq:primal-sos} \cite{lasserre_jointmarginal_2010}.

\subsubsection{Cluster basis hierarchy}\label{section:ssos:variations:cluster-basis}

In this section,
we propose a cluster basis hierarchy, wherein we utilize possible
spatial organization of the problem to sparsify
the problem and reduce the size of the SDP that must be solved \cite{vandenberghe_chordal_2017,chen_convex_2023}.
The cluster basis is a physically motivated prior often used in statistical and condensed matter physics, where we assume that our degrees of freedom can be arrayed in space, with locally close variables interacting strongly (kept in the model) and globally separated variables interacting weakly (ignored).
Moreover, one may also keep only the terms with interactions between a small number of degrees of freedom,
such as considering only pairwise or triplet interactions between particles.

In the cluster basis hierarchy, a given level of the hierarchy is defined both by the maximum degree of a variable $t$ and the desired body order $b$.
Body order denotes the maximum number of interacting
variables in a given monomial term, e.g. $x_i^a x_j^b x_k^c$ would have body order 3 and total degree $a+b+c$.
The basis function $m_{b,t}$ consists of terms $x^{\alpha}$ with $\alpha$ a multi-index,
$||\alpha||_0 \leq b$ (at most $b$ interacting variables can occur in a single term),
and $||\alpha||_{\infty} \leq t$ (each variable can have
up to degree $t$.
The maximum degree of the basis function $m_{b,t}$ is then $s = b t$.
If we are to compare $m_{b, t}$ from the cluster basis hierarchy with $m_s$ from the Lasserre hierarchy, we find
that even when $bt = s$ we still have strictly fewer terms, e.g. in the case where $b=2, t=2, s=4$ we have $m_{s}$ containing terms of the form $x_i^4$ but $m_{b,t}$ only has degree-4 terms of the $x_i^2 x_j^2$.
For further details, see discussion in \cref{appendix:snl:cluster-basis-expanded}.

\subsection{Convergence of S-SOS} \label{subsection:convergence}

As we increase the degree $s$ (either $s$ in the Lasserre hierarchy or $b, t$ in the cluster basis hierarchy) we would expect the SDP objective values $p^*_{2s}$ (\cref{eq:primal-sos-deg}) to converge to
the optimal value $p^*$ and the lower bounding function $c^*_{2s}(\omega)$ to converge to the tightest lower bound $c^*(\omega) = \inf_x f(x, \omega)$.
In this paper we refer to $p^*_{2s} \to p^*$ and $d^*_{2s} \to d^*$ interchangeably as strong duality occurs in practice
despite being difficult to formally verify (\cref{appendix:strong-duality}).
This convergence is a common feature of SOS hierarchies.
In this section we show that using polynomial $c^*_{2s}(\omega)$ to approximate $c^*(\omega)$ still allows for asymptotic convergence
in $L^1$ as $s \to \infty$. We further show how this can be improved with other choices of approximating
function classes beyond polynomial $c(\omega)$.
We specialize to the particular case of trigonometric polynomials $f(x, \omega), c(\omega)$ on $X = [0, 1]^n$ and
compact $\Omega \subset \mathbb{R}^d$ and prove asymptotic convergence of the degree-$2s$ S-SOS hierarchy as $s \to \infty$.

\subsubsection{$\ln s / s$ convergence using
a polynomial approximation to $c^*
(\omega)$}\label{section:ssos:convergence:log-s-poly}

We would like to bound the gap between the optimal lower bound $c^*(\omega) = \inf_{x\in X}f(x,\omega)$ and
the lower bound $c^*_{2s}(\omega)$ resulting from solving the degree-$2s$ primal S-SOS SDP, i.e.
\begin{align}
0 \leq c^*(\omega) - c^*_{2s}(\omega) \leq \varepsilon(f,s) \; \forall \; \omega \in \Omega. \label{eq:gap-fn}
\end{align}

To that end, we need to understand the regularity of $c^*$. Without further assumptions,
we may assume $c^*$ to be Lipschitz continuous, per \cref{proposition:lipschitz-continuity}.

With \cref{eq:gap-fn} we may then integrate
$$ 0 \leq \int_{\Omega} \inf_x f(x,\omega)-c^*_{2s}(\omega) \mathrm{d}\nu(\omega) \leq |\Omega| \varepsilon(f,s) $$ where we control $\varepsilon$ in terms of the degree $s$.
If we can drive $\epsilon \to 0$ as $s \to \infty$ then we are done.

\begin{proposition}[Theorem 2.1 in \cite{Clarke}]\label{proposition:lipschitz-continuity}
    Let $g:X\times Y\to \mathbb R$ be polynomial. Then $y\mapsto \inf_{x\in X}g(x,y)$ is Lipschitz continuous.
\end{proposition}

\begin{named}{Asymptotic convergence of S-SOS} \label{thm:convergence-abridged}
    Let $f: [0, 1]^n \times \Omega \to \mathbb{R}$ be a trigonometric polynomial of degree $2r$, $c^*(\omega) = \inf_x f(x, \omega)$ the
    optimal lower bound as a function of $\omega$,
    and $\nu$ any probability measure on compact $\Omega \subset \mathbb{R}^d$.
    Let $s$ refer to the degree of the basis in both $x, \omega$ terms and the degree of the lower-bounding polynomial $c(\omega)$, i.e.
    $m_s([x, \omega]): \mathbb{R}^{n+d} \to \mathbb{R}^{a(n,d,s)}$ is the full basis function of terms $[x, \omega]^{\alpha}$ with $||\alpha||_1 \leq s$ and $c(\omega)$ only has terms $\omega^{\alpha}$ with $||\alpha||_1 \leq s$.
    
    Let $p^*_{2s}$ be the solution to the following S-SOS SDP (c.f. \cref{eq:primal-sos-deg}) with $m_s(x, \omega)$ a spanning basis of
    trigonometric monomials with degree $\leq s$:
    \begin{align*}
        p^*_{2s} &= \sup_{c \in \mathcal{P}^{2s}(\Omega), W \succcurlyeq 0} \int c(\omega) \mathrm{d}\nu(\omega) \\
                       & \quad \text{s.t.} \quad f(x, \omega) - c(\omega) = m_s(x, \omega)^T W m_s(x, \omega) 
    \end{align*}
    
    Then there is a constant $C>0$ depending only on $||f-\bar{f}||_F, ||c^* - \bar{c}^*||, r, \Omega, n, d$ such that the following holds:
    $$ \int_{\Omega} \left[ c^*(\omega) - c^*_{2s}(\omega) \right] \mathrm{d}\nu(\omega) \leq C \frac{\ln s}{s}$$
    where $\bar{f}$ denotes the average value of the function $f$ over $[0, 1]^n$, i.e. $\bar{f} = \int_{[0,1]^n} f(x) dx$ and
    $||f(x)||_F = \sum_{\hat{x}} |\hat{f}(\hat{x})|$ denotes the norm of the Fourier coefficients.
    Thus we have asymptotic convergence of the S-SOS SDP hierarchy to the optimal value $p^*$ of
    \cref{eq:primal-sos} as we send $s \to \infty$.
\end{named}

\begin{proof}
    The following is an outline of the proof.
    For complete details, including the full theorem and proof, please see \cref{appendix:convergence-proof}.

    We define a trigonometric polynomial (t.p.)
    $c^*_a(\omega)$ of degree $s_c$
    that approximates
    the lower-bounding function 
    such that $c^*(\omega) = \inf_x f(x, \omega) \geq c^*_a(\omega)$.
    The error integral breaks apart into two terms, one bounding the approximation error
    between $c^*(\omega)$ and $c^*_a(\omega)$, and the other bounding the error between
    the approximate lower-bounding t.p. $c^*_a(\omega)$ and the SOS lower-bounding t.p. $c^*_{2s}(\omega)$.

    We then follow the proofs of \cite{fang_sum--squares_2021,bach_exponential_2023,slot_sum--squares_2023} wherein
    we define an invertible linear operator $T$ that constructs a SOS function out of a non-negative function, and show
    that such an operator exists for sufficiently large $s$.
    The core modification is to the operator $T$ which is defined as an integral operator over two kernels $q_x(x), q_{\omega}(\omega)$, i.e.
    $$ T h(x, \omega) = \int_{X \times \Omega} |q_x(x - \bar{x})|^2 |q_{\omega}(\omega - \bar{\omega})|^2 h(\bar{x}, \bar{\omega}) \mathrm{d}\bar{x} \mathrm{d}\bar{\omega} $$
\end{proof}

\subsubsection{$1/s^2$ convergence using
a piecewise-constant approximation to $c^*(\omega)$}\label{section:ssos:convergence:inv-s-piecewise}

Prior work \cite{bach_exponential_2023} achieves $1/s^2$ convergence for the regular SOS hierarchy without further assumptions.
In the previous section, we could only achieve $\ln s / s$ due to the need to first approximate the tightest lower-bounding function
$c^*(\omega)$ with a polynomial approximation, which converges at a slower rate.
To accelerate the convergence rate, we want to control the regularity of $c^*(\omega)$.
We can achieve $1/s^2$ by
approximating the $c^*(\omega)$ pointwise instead of using a smooth parameterized polynomial.
By constructing a domain decomposition of $\Omega$ and finding a SOS approximation
in $x$ for each domain, we can stitch these together to build a piecewise-constant approximation to the lower-bounding function $c^*$.

In the one-dimensional case $\Omega\subset\mathbb R$ (full proof in \cref{appendix:convergence-proof}) we achieve the following:
\begin{proposition}\label{prop:piecewise-const-convergence}
    Let $\Omega\subset\mathbb R$ be a compact interval and $f$ be a trigonometric polynomial of degree $2r$. Let $\{\omega_i\}$ be equi-distant grid points in $\Omega$ and $s_p$ the number of such points. Denote by $c^*_s(\omega_i)$ the best SOS approximation of degree $s$ of $x \mapsto f(x, \omega_i)$ and define
    $$ c^*_s = \sum_{i=1}^{s_p} c^*_s(\omega_i)1_{[\omega_i,\omega_{i+1}]}.$$ Then we have for some constant $C'$
    depending only on $\max_{\omega_i} ||f(\omega_i, \cdot) - \bar{f}(\omega_i, \cdot)||_F, r, n, \Omega, s_p$:
    \begin{align*}
        \int_\Omega c^*(\omega) - c^*_s(\omega) \mathrm{d}\omega
        \leq \max_{\omega_i}\|f(\omega_i,\cdot)-\bar f(\omega_i,\cdot)\|_F \left[1-\left(1-\frac{6r^2}{s^2}\right)^{-n}\right] |\Omega|+ \frac C{s_p} \leq C' \frac{1}{s^2}
    \end{align*}
\end{proposition}

\section{Numerical experiments}\label{section:numerics}

We present two numerical studies of S-SOS demonstrating
its use in applications. The first study (\cref{section:numerics:simple-univariate}) numerically tests
how the optimal values of the SDP \cref{eq:primal-sos-deg} $p^*_{2s}$ converge to $p^*$ of the original primal \cref{eq:primal-sos} as we increase the degree.
The second study (\cref{section:numerics:snl}) evaluates the performance of S-SOS for solution extraction
and uncertainty quantification
in various sensor network localization problems.

\subsection{Simple quadratic SOS function}\label{section:numerics:simple-univariate}

As a simple illustration of S-SOS, we
test it on the SOS function
\begin{align}\label{eq:simple-quadratic}
f(x, \omega) &= (x - \omega)^2 - (\omega x)^2
\end{align}
with $x \in \mathbb{R}, \omega \in \mathbb{R}$. The lower bound $c^*(\omega) = \inf_x f(x, \omega)$ can be computed analytically as $c^*(\omega) = \omega^4 / (1 + \omega^2)$.
Assuming $\omega \sim \text{Uniform}(-1, 1)$,
we get that the objective value for the ``tightest lower-bounding''
primal problem \cref{eq:primal-sos} is $p^* = \int_{-1}^1 \frac{\omega^4}{2 (1 + \omega^2)} d\omega = \frac{\pi}{4} - \frac{2}{3} \approx 0.1187$.
For further details, see \cref{appendix:simple-potential}.

We are interested in studying the quantitative convergence of the S-SOS hierarchy numerically.
The idea is to solve the primal (dual) degree-$2s$ SDP to find
the tightest polynomial lower bound (the minimizing probability distribution)
for varying degrees $s$. As $s$ gets larger, the basis function
$m_s(x)$ gets larger and the objective value of the SDP \cref{eq:primal-sos-deg} $p^*_{2s}$
should converge to the theoretical optimal value $p^*$.

\begin{figure}
\begin{center}
    \includegraphics[width=0.90\linewidth]{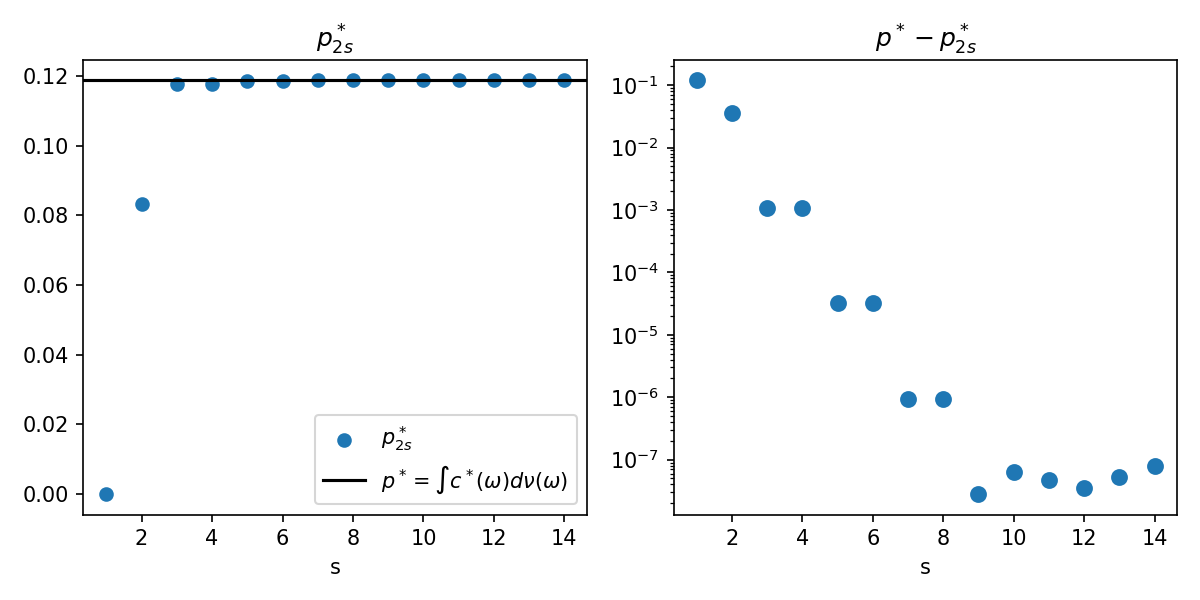}
\end{center}
\vspace{-0.3cm}
\caption{
Comparison between the objective value $p^*_{2s}$ from solving the degree-$2s$ S-SOS SDP and the objective value $p^*$ resulting from the best-possible lower bound $c^*(\omega)$
for noise drawn as $\omega \sim \text{Uniform}(-1, 1)$.
$p^* = \int c^*(\omega) d\nu(\omega) = \frac{\pi}{4} - \frac{2}{3} \approx 0.1187$ is plotted as the line in black and the
$p^*_{2s}$ values are shown as blue dots (left)
with the gap between the values $p^* - p^*_{2s}$ (right).
}
\label{fig:obj-value-vs-degree-two-plots}
\end{figure}

In \cref{fig:obj-value-vs-degree-two-plots} we see 
very good agreement between $p^*$ and $p^*_{2s}$
with exponential convergence as $s$ increases.
This is much faster than the rate we found in \cref{section:ssos:convergence:log-s-poly},
but agrees with the exponential convergence results
from \cite{bach_exponential_2023} achieved
with local optimality assumptions.
Due to the simplicity of \eqref{eq:simple-quadratic}, it's not surprising that
we see much faster convergence.
In fact, for most typical functions,
we might expect convergence much faster than the worst-case rate.
The tapering-off of the convergence rate is likely attributed to
the numerical tolerance used in our solver (CVXPY/MOSEK),
as we observed that increasing the tolerance shifts the best-achieved gap higher.

\subsection{Sensor network localization}\label{section:numerics:snl}

Sensor network localization (SNL) is a common testbed for global optimization
and SDP solvers due to the high sensitivity and ill-conditioning of the problem.
In SNL, one seeks to recover the positions of $N$ sensors
$X \in \mathbb{R}^{N \times \ell}$
positioned in $\mathbb{R}^{\ell}$ given a set of noisy observations of pairwise 
distances $d_{ij} = ||x_i - x_j||$ between the sensors \cite{nie_sum_2009,so_theory_2007}.
To have a unique global minimum and remove symmetries,
sensor-anchor distance observations are often added, where several sensors are anchored
at known locations in the space.
This can improve the conditioning of the problem,
making it ``easier'' in some sense.

\subsubsection{Definitions}

We define a SNL \emph{problem instance}
with $X \in [-1, 1]^{N \times \ell}$ as
the ground-truth positions for $\mathcal{S} = \{ 1, 2, \ldots, N \}$ sensors, $A \in [-1, 1]^{K \times \ell}$ as the ground-truth positions
for $\mathcal{A} = \{ 1, 2, \ldots, K \}$ anchors,
$\mathcal{D}_{ss}(r) = \{ d_{ij} = ||x_i - x_j|| : i, j \in \mathcal{S} \; \text{and} \; d_{ij} \leq r \}$ as the set of observed
sensor-sensor distances and $\mathcal{D}_{sa}(r) = \{ d_{ik} = ||x_i - a_k|| : i \in \mathcal{S}, k \in \mathcal{A} \; \text{and} \; d_{ik} \leq r \}$ as the set of observed sensor-anchor distances, both of which depend on some sensing radius $r$.

Writing $x_i, a_k \in [-1, 1]^{\ell}$ as the unknown positions of the $i$-th sensor and the $k$-th anchor, we can write the
potential function to be minimized as a polynomial:
\begin{align}\label{eq:snl-soft} 
f(x, \omega; X, A, r) &= \underbrace{
    \sum_{d_{ij} \in \mathcal{D}_{ss}(r)} ( ||x_i - x_j||_2^2 - d_{ij}(\omega)^2 )^2
}_{\text{sensor-sensor interactions}}
+ \underbrace{
    \sum_{d_{ik} \in \mathcal{D}_{sa}(r)} ( ||x_i - a_k||_2^2 - d_{ik}(\omega)^2 )^2
}_{\text{sensor-anchor interactions}}
\end{align}
The observed sensor-sensor and sensor-anchor distances $d_{ij}(\omega), d_{ik}(\omega)$ can be perturbed arbitrarily, but in this paper we focus on linear uniform noise, i.e. for a subset of observed
distances we have
$d_{ij,k}(\omega) = d_{ij}^* + \epsilon \omega_k$ with
$\omega_k \sim \text{Uniform}(-1, 1)$.
Other noise types may be explored, including those including outliers, which may be a better fit for robust methods (\cref{appendix:snl:noise-types}).

\cref{eq:snl-soft} contains soft penalty terms for sensor-sensor terms and sensor-anchor terms.
We can see that this is a degree-4 polynomial in the standard monomial basis elements,
and a global minimum of this function is achieved at $f(X, \mathbf{0}^d; X, A, r) = 0$ (where the distances
have not been perturbed by any noise).
In general for non-zero $\omega$ (measuring distances under noise perturbations)
we expect the function minimum to be $> 0$, as
there may not exist a configuration of sensors $\hat{X}$ that is consistent with the observed
noisy distances.

We can also support equality constraints in our solution, in particular hard equality constraints on the positions of certain sensors
relative to known anchors. This corresponds to removing all sensor-anchor soft penalty terms from the function and instead selecting
$N_H < N$ sensors at random to exactly fix in known positions via equality constraints in the SDP. The SDP is still large but the effective number of variable sensors has been reduced to
$N' = N - N_H$.

A given SNL \emph{problem type} is specified by a spatial dimension $\ell$, $N$ sensors, $K$ anchors, a sensing radius $r \in (0, 2 \sqrt{\ell})$,
a noise type (linear), and anchor type (soft penalty or hard equality). Once these are specified, we generate a random \emph{problem instance}
by sampling $X \sim \text{Uniform}(-1, 1)^n, A \sim \text{Uniform}(-1, 1)^d$. The potential $f(x, \omega)$ for a given instance
is formed (either with sensor-anchor terms or not, with terms kept based on some sensing radius $r$, and noise variables appropriately added).

The number of anchors is chosen to be as few as possible so as to still enable exact localization,
i.e. $K = \ell + 1$ anchors for a SNL problem in $\ell$ spatial dimensions. The SDPs are formulated with the help of SymPy \cite{sympy} and
solved using CVXPY \cite{diamond2016cvxpy,agrawal2018rewriting} and Mosek \cite{mosek} on a server with two Intel Xeon 6130 Gold processors (32 physical cores total) and 256GB of RAM.
For an expanded discussion and further details, see \cref{appendix:snl}.

\subsubsection{Evaluation metrics}

The accuracy of the recovered solution is of primary interest, i.e. our primary evaluation metric should be the distance between our extracted sensor positions $x$ and the
ground-truth sensor positions $X$, i.e. $\text{dist}(x, X)$.
Because the S-SOS hierarchy recovers estimates of the sensor positions $\mathbb{E}[x_i]$ along with uncertainty estimates $\text{Var}[x_i]$,
we would like to measure the distance between our ground-truth positions $X$ to our estimated distribution $p(x) = \mathcal{N}(\mathbb{E}[x], \text{Var}[x])$.
The Mahalanobis distance $\delta_M$ (\cref{eq:m-distance}) is a modified distance metric that accounts for the uncertainty \cite{mahalanobis_generalized_1936}. We use this as our primary metric for sensor recovery accuracy.

\begin{align}\label{eq:m-distance}
    \delta_M(X, \mathcal{N}(\mu, \Sigma)) &:= \sqrt{(X-\mu)^T \Sigma^{-1} (X-\mu)} 
\end{align}

As our baseline method, for each problem instance we apply a basic Monte Carlo method detailed in \cref{alg:mcpo} (\cref{appendix:snl:algorithms}) where we sample $\omega \sim \nu(\omega)$,
use a local optimization solver to find $x^*(\omega) = \inf_x f(x, \omega)$, and use this to estimate
$\mathbb{E}_{\omega \sim \nu}[x], \text{Var}_{\omega \sim \nu}[x]$.
Note that though this non-SOS method achieves some estimate of the dual SDP objective $\int f(x, \omega) d\mu(x, \omega)$, it is not guaranteed to be a lower bound.

\subsubsection{Results}

\textbf{Recovery accuracy.}
In Table \ref{results:table-ssos-mcpo} we see a comparison of the S-SOS method and the MCPO baseline.
Each row corresponds to one SNL problem type, i.e. we fix the physical dimension $\ell$, the number of anchors $K=\ell+1$, and select the sensing radius $r$ 
and the noise scale $\epsilon$. We then generate $L=20$ random instances of each problem type, corresponding to a random realization of the ground-truth sensor and anchor
configurations $X \in [-1, 1]^{N \times \ell}, A \in [-1, 1]^{K \times \ell}$, producing a $f(x, \omega)$ that we then solve the SDP for (in the case of S-SOS) or do
pointwise optimizations for (in the case of MCPO).
Each method outputs estimates for the sensor positions and uncertainty around it as a $\mathcal{N}(\mathbb{E}[x], \text{Cov}[x])$, which we then compute $\delta_M$ %
for (see \cref{eq:m-distance}), treating each dimension as independent of each other (i.e. $X$ as a flat vector).
Each instance solve gives us one observation of $\delta_M$ or each method, and we report the median and the $\pm 1\sigma_{34\%}$ values over the $L=20$ instances we generate.

\begin{table*}[!htb]
\centering
\caption{Comparison of S-SOS and MCPO solution extraction accuracy.
We present the Mahalanobis distance $\delta_M$ (\cref{eq:m-distance}) of the
the true sensor positions $X^*$ to the extracted
distribution $\mathcal{N}(\mathbb{E}[x], \text{Var}[x])$ over solutions recovered from S-SOS for varying SNL problem types.
$\ell$ is the spatial dimension, $r$ is the sensing radius used to cutoff terms in the potential $f(x, \omega)$,
$\epsilon$ is the noise scale, $N_H$ is the number of hard equality constraints used (sensors fixed
at known locations), $N_C$ is the number of clusters used (see \cref{appendix:snl:cluster-basis-expanded}), and $N$ is the number of sensors used.
Each SNL problem instance has $K = \ell + 1$ anchors used in the potential (if $N_H=0$).
The MCPO values are estimated with $T=50$ Monte Carlo iterates.
Each entry is $\hat{\mu} \pm \hat{\sigma}$ where $\hat{\mu}$ is the median and robust standard-deviation ($\sigma_{34\%}$)
estimated over 20 runs of the same problem type with varying random initializations of the sensor positions.
The entries with the lowest median $\delta_M$ are bolded.
We also compare the number of elements in the full basis $a_f$, the cluster basis $a_c$, and the
reduction multiple when using the cluster basis $a_f/a_c$.
When passing to the cluster basis, $a_f/a_c$ is how much the semidefinite matrix shrinks by.
}
\label{results:table-ssos-mcpo}
\vspace{0.1in}

\begin{tabular}{cccccc ccc cc}
\toprule
\multicolumn{6}{c}{Parameters} & \multicolumn{3}{c}{Basis comparison} & \multicolumn{2}{c}{M-distance ($\delta_M$)} \\
\cmidrule(lr){1-6} \cmidrule(lr){7-9} \cmidrule(lr){10-11}
\multicolumn{1}{r}{$\ell$} & $r$ & $\epsilon$ & $N_H$ & $N_C$ & $N$ & $a_f$ & $a_c$ & $a_f/a_c$ & S-SOS & MCPO \\
\midrule

1 & 0.5 & 0.3 & 0 & 1 & 10 & 78 & 78 & 1x & \textbf{$\boldsymbol{0.94 \pm 0.22}$} & $2.61 \pm 3.86$ \\
1 & 1.0 & 0.3 & 0 & 1 & 10 & 78 & 78 & 1x & \textbf{$\boldsymbol{0.29 \pm 0.16}$} & $1.10 \pm 0.58$ \\
1 & 1.5 & 0.3 & 0 & 1 & 10 & 78 & 78 & 1x & \textbf{$\boldsymbol{0.11 \pm 0.11}$} & $0.86 \pm 0.52$ \\
\midrule

1 & 1.5 & 0.3 & 2 & 1 & 10 & 78 & 78 & 1x & \textbf{$\boldsymbol{0.24 \pm 0.37}$} & $1.06 \pm 1.28$ \\
1 & 1.5 & 0.3 & 4 & 1 & 10 & 78 & 78 & 1x & \textbf{$\boldsymbol{0.10 \pm 0.03}$} & $0.61 \pm 0.41$ \\
1 & 1.5 & 0.3 & 6 & 1 & 10 & 78 & 78 & 1x & \textbf{$\boldsymbol{0.06 \pm 0.04}$} & $0.48 \pm 0.32$ \\
1 & 1.5 & 0.3 & 8 & 1 & 10 & 78 & 78 & 1x & \textbf{$\boldsymbol{0.04 \pm 0.02}$} & $0.31 \pm 0.17$ \\
\midrule

2 & 1.5 & 0.1 & 0 & 9 & 9 & 406 & 163 & 2.5x & \textbf{$\boldsymbol{2.86 \pm 0.94}$} & $1562.39 \pm 596.29$ \\
2 & 1.5 & 0.1 & 0 & 9 & 15 & 820 & 317 & 2.6x & \textbf{$\boldsymbol{3.25 \pm 1.19}$} & $1848.65 \pm 650.45$ \\
\bottomrule
\end{tabular}

\end{table*}

\section{Discussion}\label{section:discussion}

In this paper, we discuss the stochastic sum-of-squares (S-SOS) method to solve global polynomial optimization
in the presence of noise, prove two asymptotic convergence results for polynomial $f$ and compact $\Omega$,
and demonstrate its application to parametric polynomial minimization and uncertainty quantification
along with a new cluster basis hierarchy that enables
S-SOS to scale to larger problems.
In our experiments, we specialized to sensor
network localization and low-dimensional uniform random noise with small $n, d$.
However, it is relatively straightforward to extend
this method to support other noise types (such as Gaussian random variates without compact support, which we do in \cref{appendix:simple-potential:different-noise-distributions})
and support higher-dimensional noise with $d \gg 1$.

Scaling this method to larger problems $n \gg 1$ is an open problem for all SOS-type methods.
In this paper, we take the approach of sparsification, by making the cluster basis assumption to build up a block-sparse $W$.
We anticipate that methods that leverage
sparsity or other structure in $f$ will be promising avenues of research, as well as approximate solving methods
that avoid the explicit materialization of the matrices $W, M$.
For example, we assume that the ground-truth polynomial possesses the block-sparse structure because our SDP explicitly requires the polynomial $f(x, \omega)$ to exactly decompose into some
lower-bounding $c(\omega)$ and SOS $f_{\text{SOS}}(x, \omega)$.
Relaxing this exact-decomposition assumption and generalizing beyond polynomial $f(x, \omega), c(\omega)$ may require novel approaches and
would be an exciting area for future work.

\bibliographystyle{plainnat}
\bibliography{references}

\newpage
\appendix
\section{Appendix / supplemental material}\label{section:appendix}

\subsection{Notation} \label{appendix:notation}

Let $\mathcal{P}(X)$ and $\mathcal{P}(\Omega)$ denote the spaces of polynomials on $X \subseteq \mathbb{R}^n$
and $\Omega \subseteq \mathbb{R}^d$, respectively,
where $X$ and $\Omega$ are (not-necessarily compact) subsets of their respective ambient spaces $\mathbb{R}^n$ and $\mathbb{R}^d$.
Specifically, all polynomials of the forms below belong to their respective spaces:
$$ p(x) = \sum_{\alpha \in \mathbb{Z}_{\geq 0}} c_{\alpha} x^{\alpha} \in \mathcal{P}(X),\quad p(\omega) = \sum_{\alpha \in \mathbb{Z}_{\geq 0}} c_{\alpha} \omega^{\alpha} \in \mathcal{P}(\Omega) $$
where $x = (x_1, \ldots, x_n), \omega = (\omega_1, \ldots, \omega_d)$, $\alpha$ is a multi-index for the respective spaces, and $c_{\alpha}$
are the polynomial coefficients.

Let $\mathcal{P}^d(S)$ for some $S \in \{X, \Omega\}$ denote the subspace of $\mathcal{P}(S)$
consisting of polynomials of degree $\leq d$, i.e. polynomials where the multi-indices of the monomial terms satisfy $||\alpha||_{1} \leq d$.
$\mathcal{P}_{\text{SOS}}(X \times \Omega)$ refers to the space of polynomials on $X \times \Omega$ that can be expressible as a sum-of-squares in $x$ and $\omega$ jointly.
Additionally, $W \succcurlyeq 0$ for a matrix $W$ denotes that $W$ is symmetric positive semidefinite (PSD). Finally, $\mathbb P(\Omega)$ denotes the set of Lebesgue probability measures on $\Omega$.

\subsection{Related work}\label{appendix:related-work}

\subsubsection{Sum-of-squares theory and practice}

The theoretical justification underlying the SDP relaxations in global optimization we use here derive from the Positivstellens\"atz (positivity certificate) of \cite{putinar_positive_1993}, a representation theorem guaranteeing that strictly positive polynomials on certain sets admit sum-of-squares representations. Following this, \cite{lasserre_global_2001,lasserre_moment-sos_2018,lasserre_moment-sos_2023} developed the Moment-SOS hierarchy, describing a hierarchy of primal-dual SDPs (each having fixed degree) of increasing size that provides a monotonic non-decreasing sequence of lower bounds.

There is rich theory underlying the SOS hierarchy combining disparate results from algebraic geometry \cite{parrilo_structured_2000,lasserre_moment-sos_2018,lasserre_moment-sos_2023},
semidefinite programming \cite{nie_sum_2009,papp_sum--squares_2019},
and complexity theory \cite{de_klerk_complexity_2008,odonnell_sos_2016}.
The hierarchy exhibits finite convergence in particular cases where convexity and a strict local minimum are guaranteed
\cite{nie_optimality_2014}, otherwise
converging asymptotically \cite{bach_exponential_2023}.
In practice, the hierarchy often
does even better than these guarantees, converging exactly at $c^*_s$ for some small $s$.

The SOS hierarchy has found numerous applications in wide-ranging fields, including: reproducing certain results of perturbation theory and providing useful lower-bound certifications in quantum field theory and quantum chemistry \cite{hastings_perturbation_2022,hastings_field_2023}, providing better provable guarantees in high-dimensional statistical problems \cite{hopkins_statistical_2018,hopkins_mixture_2018}, useful applications in the theory and practice of sensor network localization \cite{nie_sum_2009,sedighi_localization_2021} and in robust and stochastic optimization \cite{bertsimas_hierarchy_2011}.

Due to the SDP relaxation, the SOS hierarchy is quite powerful. This flexibility comes at a cost, primarily in the form of computational complexity.
The SDP prominently features a PSD matrix $W \in \mathbb{R}^{a(n,d,s) \times a(n,d,s)}$ with $a(n,d,s)$ scaling as $\binom{n+d+s}{s}$ for $n$ dimensions and maximum degree $s$.
Without exploiting the structure of the polynomial, such as locality (coupled terms) or sparsity, solving the SDP using a standard interior point method becomes prohibitively expensive for moderate values of $s$ or $n$. Work attempting to improve the scalability of the core ideas underlying the SOS hierarchy and the SDP method include \cite{ahmadi_dsos_2019,papp_sum--squares_2019}.

\subsubsection{Stochastic sum-of-squares and parametric polynomial optimization}

The S-SOS hierarchy we present in this paper as a solution to parametric polynomial optimization was presented originally by \cite{lasserre_jointmarginal_2010} as a ``Joint + Marginal'' approach. That work provides the same hierarchy of semidefinite
relaxations where the sequence of optimal solutions converges to
the moment vector of a probability measure encoding all information
about the globally-optimal solutions $x^*(\omega) = \text{argmin}_{x} f(x, \omega)$ and provides a proof that the dual problem (our primal)
obtains a polynomial approximation to the optimal value function
that converges almost-uniformly to $c^*(\omega)$.

\subsubsection{Uncertainty quantification and polynomial chaos}

Once a physical system or optimization problem is characterized, sensitivity analysis and uncertainty quantification seek to quantify how randomness or uncertainty in the inputs can affect the response.
In our work, we have the parametric problem of minimizing a function $f(x, \omega)$ over $x$ where $\omega$ parameterizes the function and is drawn from some noise distribution $\nu(\omega)$.

If only function evaluations $f(x, \omega)$ are allowed and no other information is known, Monte Carlo is often applied, where one draws $\omega_k \sim \nu(\omega)$ and solves many realizations of $\inf_x f_k(x) = f(x, \omega_k)$ to approximately solve the
following stochastic program:
$$ f^* = \inf_x \mathbb{E}_{\omega \sim \nu}[f(x, \omega)]$$
Standard Monte Carlo methods are ill-suited for integrating high-dimensional functions, so this method is computationally challenging in its own right.
In addition, we have no guarantees on our result except that as we take the number of Monte Carlo iterates $T \to \infty$ we converge to some unbiased
estimate of $\mathbb{E}_{\omega \sim \nu}[f(x, \omega)]$.

Our approach to quantifying the uncertainty in optimal function value resulting from uncertainty in parameters $\omega$ is to find a deterministic lower-bounding $c^*(\omega)$ which guarantees $f(x, \omega) \geq c^*(\omega)$ no matter the realization of noise.
This is reminiscent of the polynomial chaos expansion literature, wherein a system of some stochastic variables is expanded into a deterministic function of 
those stochastic variables, usually in some orthogonal polynomial basis \cite{sudret_global_2008,najm_uncertainty_2009}.

\subsection{An example}\label{appendix:sdp-examples}

\begin{example}
Let $f(x, \omega)$ be some polynomial of degree $\leq 2s$
written in the standard monomial basis, i.e.
\begin{align*}
f(x, \omega) &= \sum_{||\alpha||_{1} \leq 2s} f_{\alpha} x^{\alpha} \\
&= \sum_{||\alpha||_{1} \leq 2s}
f_{(\alpha_1, \ldots, \alpha_{n + d})}
\prod_{i=1}^n x_i^{\alpha_i}
\prod_{i=1}^d \omega_i^{\alpha_{n+i}}
\end{align*}
Let $m_s(x, \omega) \in \mathbb{R}^{a(n,d,s)}$ be the 
basis vector representing the full set of monomials in
$x, \omega$ of degree $\leq s$ with
$a(n,d,s) = \binom{n+d+s}{s}$.

For all $\alpha \in \mathbb{Z}_{\geq 0}^{n+d}$ with $||\alpha||_{1} \leq 2s$
and $\alpha_k = 0$ for all $k \in \{1, \ldots, n\}$ (i.e.
monomial terms containing only $\omega_1, \ldots, \omega_d$)
we must have:
$$ \int_{X \times \Omega} \omega^{\alpha} \mathrm{d}\mu(x, \omega) - \int_{\Omega} \omega^{\alpha} \mathrm{d}\nu(\omega) = 0$$

Explicitly, for $\mu$ to be a valid probability distribution we must have:
$$ \int_{X \times \Omega} \mathrm{d}\mu(x, \omega) - 1 = M_{0, 0} - 1 = y_{(1, 0, \ldots)} - 1 = 0 $$
Suppose $\Omega = [-1, 1], \omega \sim \text{Uniform}(-1, 1)$ so that $d=1$, $\nu(\omega) = 1/2$.
We require:
\begin{align*}
\int_{X \times \Omega} \omega^{\alpha} \mathrm{d}\mu(x, \omega) = \int_{[-1, 1]} \omega^\alpha \mathrm{d}\nu(\omega) = 
    \begin{cases}
    1 & \alpha = 0 \\
    0 & \alpha = 1 \\
    \frac{1}{3} & \alpha = 2 \\
    0 & \alpha = 3 \\
    \frac{1}{5} & \alpha = 4 \\
    \end{cases}
\end{align*}
\qed
\end{example}

\subsection{Strong duality}\label{appendix:strong-duality}

To guarantee strong duality theoretically, we need a strictly feasible point in the interior (Slater's condition).
For us, this is a consequence of Putinar's Positivstellensatz,
if $f(x, \omega)$ admits a decomposition as $f(x, \omega) = c(\omega) + g(x, \omega)$
where $g(x, \omega) > 0$ (i.e. is strictly positive), we have strong duality, i.e. $p^* = d^*$ and $p^*_{2s} = d^*_{2s}$ \cite{lasserre_global_2001,schmudgen_moment_2017}.
However, it is difficult to verify the conditions analytically. In practice, strong duality is observed in most cases,
so in this paper we refer to solving the primal and dual interchangeably, as $p^*_{2s} = d^*_{2s}$ in all cases we encounter
where a SDP solver returns a feasible point.

\subsection{Proofs}

\subsubsection{Primal-dual relationship of S-SOS}\label{appendix:formal-dual}

\paragraph{Regular SOS}
Global polynomial optimization can be framed as the following lower-bound maximization problem
where we need to check global non-negativity:
\begin{align}\label{eq:appendix-reg-primal}
    \sup_{c\in\mathbb R} \quad & c \\
    \text{s.t.} \quad & f(x) - c \geq 0 \quad \forall x \notag
\end{align}
When we take the SOS relaxation of the non-negativity constraint in the primal,
we now arrive at the SOS primal problem, where we require $f(x) - c$ to be SOS
which guarantees non-negativity but is a stronger condition than necessary:
\begin{align}\label{eq:appendix-reg-primal-sos}
    \sup_{c\in\mathbb R} \quad & c \\
    \text{s.t.} \quad & f(x) - c \in \mathcal P_{\rm SOS}(X). \notag
\end{align}

The dual to \cref{eq:appendix-reg-primal} is the following moment-minimization problem:
\begin{align}\label{eq:appendix-reg-dual}
    \inf_{\mu\in\mathbb P(X)} \quad & \int f(x) \mathrm{d}\mu(x) \\
    \text{with} \quad & \int \mathrm{d}\mu(x) = 1. \notag
\end{align}
Taking some spanning basis $m_s(x): \mathbb{R}^n \to \mathbb{R}^{a(n,s)}$ of monomials up to degree $s$, we have the moment matrix $M \in \mathbb{R}^{a(n,s) \times a(n,s)}$:
$$ M_{i, j} = \int m_i(x) m_j(x) d\mu(x) = y_{\alpha} $$
where we introduce a moment vector $y$ whose elements correspond to the unique moments of the matrix $M$.
Then we may write the degree-$2s$ moment-minimization problem, which is now in a solvable numerical form:
\begin{align}\label{eq:appendix-reg-dual}
    \inf_{y} \quad & \sum_{\alpha} f_{\alpha} y_{\alpha} \\
    \text{with} \quad & M(y)_{1, 1} = 1 \notag \\
    & M(y) \succcurlyeq 0 \notag
\end{align}
where we write $M(y)$ as the matrix formed by placing the moments from $y$ into their appropriate places and we set the first element of $m_s(x)$ to be 1, hence $M_{1, 1} = \int d\nu(x) = 1$ is simply the normalization constraint. For further reading, see \cite{nie_sum_2009,lasserre_global_2001}.

\paragraph{Stochastic SOS}
Now let us lift this problem into the stochastic setting with parameters $\omega$ sampled from a given distribution $\nu$, i.e. replacing $x \to (x, \omega)$.
We need to make some choice for the objective. The expectation of the lower bound under $\nu(\omega)$
is a reasonable choice, i.e.
$$ \int_{\Omega} c(\omega) d\nu(\omega) $$
but we could also make other choices, such as ones that encourage more robust lower bounds.
In this paper however, we formulate the primal S-SOS as below (same as \cref{eq:primal-sos}):
\begin{align}\label{eq:appendix-prim-stoch}
    p^* = \sup_{c \in L^1(\Omega)}        \quad & \int c(\omega) \mathrm{d}\nu(\omega) \\
                   \text{s.t.} \quad & f(x, \omega) - c(\omega) \geq 0 \notag
\end{align}

Note that if the ansatz space for the function $c(\omega)$ is general enough, the maximization of the curve $c$ is equivalent to a pointwise maximization, i.e. we recover the best approximation for almost all $\omega.$ Then the dual problem has a very similar form to the non-stochastic case.

\begin{theorem}\label{thm:dual_form_stoch}
The dual to \cref{eq:appendix-prim-stoch} is the following moment minimization where $\mu(x, \omega)$ is a probability measure on $X \times \Omega$:
  \begin{align*}
    \inf_{\mu\in \mathbb P(X\times\Omega)}& \quad  \int f(x, \omega) \mathrm{d}\mu(x, \omega) \\
    \text{with} \quad & \int_{X \times \Omega} \omega^\alpha \mathrm{d}\mu(x,\omega) = \int_\Omega \omega^\alpha \mathrm{d}\nu(\omega)\quad \text{for all } \alpha\in\mathbb N^d.
\end{align*}  
\end{theorem}

\begin{remark}
    Notice, that the condition $\int_{X \times \Omega} \omega^\alpha \mathrm{d}\mu(x,\omega) = \int_\Omega \omega^\alpha \mathrm{d}\nu(\omega)$ implies that the first marginal of $\mu$ is the noise distribution $\nu$. Let $\mu_\omega$ denote the disintegration of $\mu$ with respect to $\nu$, \cite{AmbrosioGigliSavare}. Then the moment matching condition is equivalent to $\mu_\omega(X) = 1$ for almost all $\omega$ and $\mu$ being a Young measure w.r.t. $\nu$. The idea is that $\mu_{\omega}(x)$ is a minimizing density for every single configuration of $\omega$.
\end{remark}

\begin{proof}
We use $\mathcal{P}_{\geq 0}(X \times \Omega)$ to denote the space of non-negative polynomials on
$X \times \Omega$.
Given measure $\nu$ on $\Omega$ and polynomial function $p:X\times \Omega\to\mathbb R$ consider
\begin{align*}
    &\sup_{\substack{\gamma\in L^1(\Omega,\nu)\\ q\in \mathcal P_{\geq 0}(X\times\Omega).}} \int_{\Omega} \gamma(\omega)\mathrm{d}\nu(\omega)\\
    \mathrm{s.t}&\quad  p(x,\omega)-\gamma(\omega) = q(x,\omega)
\end{align*}
This is equivalent to
\begin{align*}
-\inf_{\substack{\gamma\in L^1(\Omega,\mu)\\q\in \mathcal P_{\geq 0}(X\times\Omega)}} f(\gamma,q)+g(\gamma,q)
\end{align*}
with $$f(\gamma,q) = -\int_{\Omega}\gamma(\omega)\mathrm{d}\nu(\omega)$$ and $$g(\gamma,q) = -\chi_{\{f-\gamma-q=0\}}=\begin{cases}
    0\text{ if }f-\gamma-q=0\\
    -\infty \text{ else}
\end{cases},$$ 
i.e. $g$ is the characteristic function enforcing non-negativity.

Denote by $h^*$ the Legendre dual, i.e. $$h^*(y) = \sup_{x}\langle x,y\rangle -h(x).$$
Then by Rockafellar duality, \cite{Ekeland,Rockafellar}, and noting that signed Borel measures $\mathcal B$ are the dual to continuous functions, the dual problem reads
\begin{align*}
    \sup_{\Gamma\in L^\infty(\Omega,\mu),\mu\in \mathcal B} -f^*(\Gamma,\mu)-g^*(-(\Gamma,\nu))
\end{align*}
and we would have 
$$\sup_{\Gamma\in L^\infty(\Omega,\mu),\mu\in \mathcal B} -f^*(\Gamma,\mu)-g^*(-(\Gamma,\mu)) = -\inf_{\substack{\gamma\in L^1(\Omega,\mu)\\q\in \mathcal P_{\geq 0}(X\times\Omega)}} f(\gamma,q)+g(\gamma,q).$$

The Legendre duals of $f$ and $g$ can be explicitly calculated as
 $$f^*(\Gamma,\mu) = \begin{cases}0\text{ if } \Gamma =-1 \text{ and } \mu\leq0\\\infty \text{ else}\end{cases}$$
 and
 $$g^*(\Gamma,\mu) = \begin{cases}
     \displaystyle \int_{\Omega\times X} f(x,\omega)\mathrm{d}\mu(\omega,x) & \text{ if } f-\gamma\in \mathcal P_{\geq 0}(X\times\Omega) \text{ and } \Gamma(\omega) = \mu_\omega(X)\\
     \infty & \text{ else}
 \end{cases}$$
since 
\begin{align*}
    f^*(\Gamma,\mu) &= \sup_{\gamma,q} \left(\int_\Omega \gamma(\omega)\Gamma(\omega)\mathrm{d}\nu(\omega) + \int_{\Omega\times X}q(x,\omega)\mathrm{d}\mu(x,\omega)
-f(\gamma,q)\right)\\
 &=\sup_{\gamma,q} \int_{\Omega} \gamma(\omega)(\Gamma(\omega)+1)\mathrm{d}\nu(\omega)+\int_{\Omega\times X}q(x,\omega)\mathrm{d}\mu(x,\omega)\\
 &= \begin{cases}
     \displaystyle 0 \quad \text{if } \Gamma=-1 \text{ and } \mu\leq 0\\
     \infty \quad \text{else}
 \end{cases}\end{align*}
 and 
 \begin{align*}
     g^*(\Gamma,\mu) & = \sup_{\gamma,q} \int_\Omega\gamma(\omega)\Gamma(\omega)\mathrm{d}\nu(\omega)+\int_{\Omega\times X}q(x,\omega)\mathrm{d}\mu(\omega,x) + \chi_{\{f-\gamma-q=0\}}\\
     &= \begin{cases}
         \displaystyle
         \sup_{\gamma} \int_\Omega\gamma(\omega)\Gamma(\omega)\mathrm{d}\nu(\omega)+\int_{\Omega\times X}(f(x,\omega)-\gamma(\omega))\mathrm{d}\mu(\omega,x)\quad &\text{if }f-\gamma=\in \mathcal P_{\geq 0}(X\times\Omega)\\\infty\quad &\text{else}
     \end{cases}\\
     &= \begin{cases}
         \displaystyle
         \sup_{\gamma} \int_\Omega\gamma(\omega)(\Gamma(\omega)-\mu_\omega(X))\mathrm{d}\nu(\omega)+\int_{\Omega\times X}(f(x,\omega)\mathrm{d}\mu(\omega,x)\quad &\text{if }f-\gamma\in \mathcal P_{\geq 0}(X\times\Omega)\\\infty\quad &\text{else}
     \end{cases}\\
          &= \begin{cases}
          \displaystyle
          \int_{\Omega\times X}(f(x,\omega)\mathrm{d}\mu(\omega,x)\quad &\text{if }f-\gamma\in \mathcal P_{\geq 0}(X\times\Omega) \text{ and }\Gamma(\omega) = \mu_\omega(X)   \\\infty\quad &\text{else}
     \end{cases}
 \end{align*}
 Altogether, we get 
 \begin{align*}
     -f^*(\Gamma,\mu)-g^*(-\Gamma,-\mu) &  = \begin{cases}
         \displaystyle
         \int_{\Omega\times X}f(x,\omega)\mathrm{d}\mu(\omega,x) \quad &\text{if }\mu_\omega(X)=1\\
         \infty\quad &\text{else.}
     \end{cases}
 \end{align*}
\end{proof}

\subsubsection{Convergence of S-SOS hierarchy}\label{appendix:convergence-proof}

\paragraph{Lemma on approximating polynomials}
\begin{lemma}\label{lem:Fourier}
    Let $\Omega$ be compact and $g: \Omega\to\mathbb R^n$ be Lipschitz continuous. Then there is a trigonometric polynomial $g_s$ of degree $s$ and a constant $C > 0$ depending only on $\Omega$ and $n$ such that $$g\geq g_s$$ and $$\|g-g_s\|_{L^2(\Omega)}\leq \frac{1+\ln(s)}{s}C\|g\|_{H^1(\Omega)}.$$
\end{lemma}

One cannot expect much more as the following example shows:
\begin{example}
    Consider $g:\mathbb R\times\mathbb R^2 \to \mathbb R$ defined by $$g(x,p,q) = (x^2+px+q)^2.$$
    Then we have for every $(p,q)\in\mathbb R^2$ that \begin{align*}
        \inf_{x\in\mathbb R} g(x,p,q) = \begin{cases}
            0 & \text{if } \frac{p^2} 4\geq q\\
            (\frac {p^2}{4} - q)^2 & \text{else.}
        \end{cases}
    \end{align*}
    Therefore, $(p,q)\mapsto \inf_{x\in\mathbb R} g(x,p,q)$ is once differentiable but not twice.
\qed
\end{example}

\paragraph{Convergence at $\ln s / s$ rate}

\begin{named}{Asymptotic convergence of S-SOS} \label{thm:convergence}
    Let $f: [0, 1]^n \times \Omega \to \mathbb{R}$ be a trigonometric polynomial of degree $2r$, $c^*(\omega) = \inf_x f(x, \omega)$ the
    optimal lower bound as a function of $\omega$,
    and $\nu$ any probability measure on compact $\Omega \subset \mathbb{R}^d$.
    Let $s = (s_x, s_\omega, s_c)$, referring separately to the degree of the basis in $x$ terms, the degree of the basis in $\omega$ terms,
    and the degree of the lower-bounding polynomial $c(\omega)$.
    
    Let $c^*_{2s}(\omega)$ be the lower bounding function obtained from the primal S-SOS SDP with $m_s(x, \omega)$ a spanning basis of
    trigonometric monomials with degree $\leq s_x$ in $x$ terms and of degree $\leq s_\omega$ in $\omega$ terms:
    \begin{align*}
        p^*_{2s} &= \sup_{c \in \mathcal{P}^{2s_c}(\Omega), W \succcurlyeq 0} \int c(\omega) \mathrm{d}\nu(\omega) \\
                       & \quad \text{s.t.} \quad f(x, \omega) - c(\omega) = m_s(x, \omega)^T W m_s(x, \omega) 
    \end{align*}

    Then there is a constant $C>0$ depending only on $\Omega,d,$ and $n$ such that for all $s_\omega,s_x\geq \max\{3r,3s_c\}$ the following holds:
    $$ \int_{\Omega} \left[ c^*(\omega) - c^*_{2s}(\omega) \right] \mathrm{d}\nu(\omega) \leq |\Omega| \epsilon(f, s)$$
    \begin{align*}
        \varepsilon(f, s)
        \leq & \|f-\bar f\|_F \left[1-\left(1-\frac{6r^2}{s_\omega^2}\right)^{-d}\left(1-\frac {6r^2}{s_x^2}\right)^{-n}\right]\\
        &\ +\|c^* - \bar c^*\|_F \left[1 - \left( 1 - \frac{6r^2}{s_\omega^2} \right)^{-d} \right] 
         + C \frac{(1+\ln(2s_c))}{2s_c}.
    \end{align*}
    
    where $\bar{f}$ denotes the average value of the function $f$ over $[0, 1]^n$, i.e. $\bar{f} = \int_{[0,1]^n} f(x) dx$ and
    $||f(x)||_F = \sum_{\hat{x}} |\hat{f}(\hat{x})|$ denotes the norm of the Fourier coefficients. 

    $\epsilon(f,s)$ bounds the expected error, giving us asymptotic convergence as $s = \min(s_x, s_{\omega}, s_c) \to \infty$.
    Note the first two terms give a $O(\frac{1}{s^2})$ convergence rate.
    However, the overall error will be dominated by the degree of $c(\omega)$ (from the third term) hence our convergence rate is $O(\frac{\ln s}{s})$.
\end{named}

\begin{proof}
    By the convergence of Fourier series \cite{jackson1930theory} we have the existence of a trigonometric polynomial $g'$ of degree $s$ with $$\|g-g'\|_{L^1(\Omega)}\leq \frac {C'}s\|g\|_{H^1(\Omega)}$$ as well as $$\|g-g'\|_\infty \leq L_g\frac {\ln(s)}s.$$
    Then we define $g_s = g'-\|g-g'\|_\infty$ and hence $g \geq g_s$. Furthermore, $$\|g-g_s\|_{L^1(\Omega)}\leq \frac{(C'+|\Omega|\ln s)} s L_g. $$
    Writing $C(\Omega) = \max \{ C', |\Omega| \}$ we have the desired form where $|\Omega|$ is the volume of $\Omega.$
\end{proof}

\begin{proof}[Proof of \cref{thm:convergence}]
    Let $\Omega\subset\mathbb R^d$ be compact and $f:\mathbb \mathbb{R}^n \times \Omega \to\mathbb R$ be a 1-periodic trigonometric polynomial (t.p.) of degree $\leq 2r$.
    We then make $\Omega$ isomorphic to $[0, 1]^d$ and hereafter consider $\Omega = [0, 1]^d$
    and $f: [0, 1]^n \times [0, 1]^d \to \mathbb{R}$.
    Let $\varepsilon>0$ and $b=\frac \varepsilon  2$. Let the best lower bound be
    $$ c^*(\omega) = \inf_{x\in X} f(x,\omega). $$

    \emph{Proof outline.}
    We split the error into two parts. First, we use the fact that there is a lower-bounding t.p. $c^*_a$ of degree $s_c$ such that $$\|c^* - c^*_a\| \leq C\frac{1+\ln s_c}{s_c}$$ and $$c^* \geq c^*_a.$$ This will provide us with a degree-$s_c$ t.p. approximation to the lower bounding function, which
    in general is only known to be Lipschitz continuous.

    Next, we show, that for any $b>0$ there is a degree-$2s$ SOS t.p. $f_{\text{SOS}}(x, \omega)$ such that
    $$f_{\text{SOS}} = f - (c^*_a - b).$$
    We write $s = (s_x, s_{\omega})$ where $s_x, s_{\omega}$
    denotes the respective max degrees in the variables $x, \omega$.
    Once we have constructed this, we can compute
    $f - f_{\text{SOS}} = c^*_a - \varepsilon$ and since we know that $f_{\text{SOS}} \geq 0$ everywhere and $c^*_a - \varepsilon$ is some degree-$s_c$ t.p.
    we have found a degree-$s_c$ lower-bounding t.p.
    The construction of this SOS t.p. adds another error term.
    If we can drive $\varepsilon \to 0$ as $\bar{s} = \min(s_x, s_{\omega}, s_c) \to \infty$ then we are done.

    \emph{Proof continued.}
    To that end, let $c^*_a:\Omega\to\mathbb R$ be the best degree-$s_c$ trigonometric approximation of $c^*$ with respect to $L^1$
    such that $$c^* \geq c_a^*.$$

    By \cite{Clarke}, we know that $c^*$ is locally Lipschitz continuous with Lipschitz constant $L_{c^*}$ and hence, by \cref{lem:Fourier} we get that there is $C(\Omega)>0$ such that $$\|c^* -c_a^* \|_{L^1\Omega)}\leq C(\Omega) \frac{1+\ln s_c}{s_c} L_{c^*}.$$

    Next we introduce $c^*_{2s}(\omega)$ which is some degree-$2s$ t.p.
    After an application of the triangle inequality and Cauchy-Schwarz on the integrated error term $\int_{\Omega} |c^* - c^*_{2s}| d\omega$ we have
    $$\int_\Omega \bigg| \inf_{x\in X} f(x,\omega)-c^*_{2s}(\omega) \bigg| \mathrm{d}\omega \leq \int_\Omega |c^*_a(\omega)-c^*_{2s}(\omega)|\mathrm{d}\omega + {|\Omega|} \|c^* -c^*_a\|_{L^2(\Omega)}$$
    $$ \int_\Omega
    \bigg| \inf_{x\in X} f(x,\omega)-c^*_{2s} (\omega) \bigg|
    \mathrm{d}\omega \leq
    \underbrace{\int_\Omega |c^*_a(\omega)-c^*_{2s}(\omega)|\mathrm{d}\omega}_{\text{gap between some SDP solution $c^*_{2s}(\omega)$ and t.p. $c^*_a(\omega)$}}
    + \underbrace{C(\Omega) \frac{1 + \ln s_c}{s_c} L_{c^*}}_{\text{approx. error of L-contin. fn.}}
    $$

    Now we want to show that for any $\varepsilon > 0$ we can construct a degree-$2s$ SOS trigonometric polynomial $f_{\text{SOS}}(x, \omega)$ such that $$f_{\text{SOS}} = f - c^*_a + b.$$
    with $b = \varepsilon/2$ and $s = (s_x, s_{\omega}) > r$.
    We can then set $f - f_{\text{SOS}} = c_a^* - b = c^*_{2s}$
    as the degree-$2s$ lower-bounding function.
    If we can drive $b = \varepsilon/2 \to 0$ as $s, s_c \to \infty$ we are done, as by construction $|c_a^* - c_{2s}^*| = b$.
    
    Observe that by assumption $f - c^*_a + b$ is a t.p. in $(x, \omega)$ where $f$ is degree-$2r$ and $c^*_a$ is degree $s_c \geq 2r$. Denote by $(f-f_*^a+b)_\omega$ its coefficients w.r.t the $\omega$ basis. Note that the coefficients are functions in $x$. 
    Following the integral operator proof methodology in \cite{bach_exponential_2023}, define the integral operator $T$ to be
    $$Th(x,\omega) = \int_{X\times \Omega} |q_\omega(\omega-\bar\omega)|^2|q_x(x-\bar x)|^2h(\bar x,\bar \omega)\mathrm{d}\bar x\mathrm{d}\bar \omega,$$
    where $q_\omega$ is a trigonometric polynomial in $\omega$ of degree $\leq s_\omega$ and $q_x$ is a trigonometric polynomial in $x$ of degree $\leq s_x$.
    The intuition is that this integral operator explicitly
    builds a SOS function of degrees $(s_x, s_{\omega})$ out
    of any non-negative function $h$ 
    by hitting it against the kernels $q_{x}, q_{\omega}$.
    
    We want to find a positive function $h:X\times\Omega\to\mathbb R$ such that $$Th = f - c_a^* + b.$$
    In frequency space, the Fourier transform turns a convolution into pointwise multiplication so we have:
    $$\widehat{Th}(\hat x,\hat \omega)=\hat q_\omega*\hat q_\omega(\hat\omega)\cdot \hat q_x*\hat q_x(\hat x) \cdot \hat h(\hat x,\hat \omega).$$
    In the Fourier domain it is easy to write down the coefficients of $\hat{h}$:
    $$
    \hat h(\hat x,\hat \omega) = \begin{cases}
    0\quad & \text{if }\|\hat x, \hat\omega\|_\infty > \max\{2r, 2s_c\} \\
    \displaystyle \frac{\hat f(\hat x, \hat \omega) -\hat c^*_a(\hat\omega)1_{\hat x=0} + b 1_{\hat x=0} 1_{\hat\omega=0}}{\hat q_\omega*\hat q_\omega(\hat \omega)\cdot \hat q_x*\hat q_x(\hat x)} & \text{otherwise}. 
    \end{cases}
    $$
    
    Computing $Th - h$ gives:
    \begin{align*}
        &f(x, \omega) - c_a^*(\omega) + b - h(x,\omega) \\
        =& \sum_{\hat\omega,\hat x} \hat f(\hat x,\hat\omega) \left( 1-\frac{1}{\hat q_\omega*\hat q_\omega(\hat \omega)\cdot \hat q_x*\hat q_x(\hat x)} \right) \exp(2i\pi \hat\omega^T \omega)\exp(2i\pi\hat x ^T x)\\
        &\ + \sum_{\hat\omega} (b 1_{\hat\omega=0} - c_a^*) \left( 1 - \frac{1}{\hat q_\omega*\hat q_\omega(\hat \omega)} \right) \exp(2i\pi \hat{\omega}^T \omega)
    \end{align*}
    and thus after requiring $\hat q_\omega * \hat q_\omega(0) = \hat q_x * \hat q_x(0) = 1$ we have:
    \begin{align*}
        &\max_{x,\omega}|f(x,\omega) - c_a^*(\omega) + b - h(x, \omega)|\\
        \leq & \|f-\bar f\|_F \max_{\hat\omega\neq 0}\max_{\hat x\neq 0} \bigg| 1-\frac{1}{\hat q_\omega*\hat q_\omega(\hat \omega)\cdot \hat q_x*\hat q_x(\hat x)} \bigg| \\
        &\ + \max_{\hat\omega\neq 0} \|c_a^* - \bar c_a^* \|_F \bigg| 1-\frac{1}{\hat q_\omega*\hat q_\omega(\hat \omega)} \bigg|.
    \end{align*}

    As a reminder, because $c^* \geq c_a^*$ everywhere we have $f - c_a \geq f - c^* \geq 0$ or $f - c_a^* + b > 0$, since $b = \varepsilon / 2 > 0$.
    Since $Th = f - c_a^* + b > 0$ and it is a SOS, we need to guarantee $h > 0$.
    
    If $\max_{x,\omega}|f(x,\omega)-f_*^a(\omega)+b-h(x,\omega)|\leq b$ then $$\max_{x,\omega}{|Th-h|} < b.$$
    Since $Th\geq b$ and $b>0$ we have
    $$h=Th+h-Th\geq Th-\|h-Th\|_\infty\geq b-b\geq 0$$
    and hence $h>0$ if we ensure $\max_{x, \omega} |Th - h| \leq b$.

    Now let us show that
    $$ \max_{x,\omega}|f(x,\omega)-c_a^*(\omega) + b -h(x,\omega)|\leq b $$
    can be ensured if $s = (s_x, s_{\omega})$ is large enough.

    Using the same kernel and bounds as in \cite{bach_exponential_2023}, we choose for $z\in\{x,\omega\}$ the triangular kernel such that $$\hat q_z(\hat z)=\left( 1-\frac {6r^2}{z^2} \right)^d_{+} \prod_{i=1}^d \left(1-\frac{|\hat z_i|}{s_{x,\omega}} \right)_{+}. $$
    Note that $(x)_+ = \max(x, 0)$.
    Then we have
    \begin{align*}
        &\max_{x}|f(x,\omega) - c_a^*(\omega)+b-h(x,\omega)|\\
        \leq & \|f-\bar f\|_F \max_{\hat\omega, \hat x} \bigg| 1-\frac{1}{\hat q_\omega*\hat q_\omega(\hat \omega) \cdot \hat q_x*\hat q_x(\hat x)} \bigg|
        +\|c_a^* - \bar c_a^* \|_F \max_{\hat \omega} \bigg| 1-\frac{1}{\hat q_\omega*\hat q_\omega(\hat \omega)} \bigg| \\
        \leq & \|f-\bar f\|_F \bigg|1-\left(1-\frac{6r^2}{s_\omega^2} \right)^{-d} \left(1-\frac {6r^2}{s_x^2}\right)^{-n} \bigg|+\|c_a^* - \bar c_a^*\|_F \bigg| 1- \left( 1-\frac{6^2}{s_\omega^2} \right)^{-d} \bigg|
    \end{align*}
    Therefore, by choosing $s_\omega$ and $s_x$ large enough such that
    $$ \|f-\bar f\|_F \bigg|1-\left(1-\frac{6r^2}{s_\omega^2} \right)^{-d} \left(1-\frac {6r^2}{s_x^2}\right)^{-n} \bigg|+\|c_a^* - \bar c_a^*\|_F \bigg| 1- \left( 1-\frac{6^2}{s_\omega^2} \right)^{-d} \bigg| \leq b = \frac{\varepsilon}{2}$$
   we have $$h \geq 0$$ and thus $Th$ is SOS. By design we have $$c^*_a - c^*_{2s}\leq b$$ 
    and thus 
    $$\int_{\Omega} |c^*_a - c^*_{2s}|\mathrm{d}\omega \leq \frac \varepsilon 2. $$

    Recalling $$ \int_\Omega
    \bigg| \inf_{x\in X} f(x,\omega)-c^*_{2s} (\omega) \bigg|
    \mathrm{d}\omega \leq
    \underbrace{\int_\Omega |c^*_a(\omega)-c^*_{2s}(\omega)|\mathrm{d}\omega}_{\text{gap between some SDP solution $c^*_{2s}(\omega)$ and t.p. $c^*_a(\omega)$}}
    + \underbrace{C(\Omega) \frac{1 + \ln s_c}{s_c} L_{c^*}}_{\text{approx. error of L-contin. fn.}}
    $$
    we can additionally choose $s_c$ large enough
    to guarantee $$ C(\Omega) \frac{1 + \ln s_c}{s_c} L_{c^*} \leq \frac{\varepsilon}{2} $$
    and then we are done.

    Setting $s_x, s_{\omega}, s_c = s$ and sending $s \to \infty$
    we have asymptotic behavior of the final error expression:
    $$
    \boxed{ \int_\Omega
    \bigg| \inf_{x\in X} f(x,\omega)-c^*_{2s} (\omega) \bigg|
    \mathrm{d}\omega \leq
    C_1 \frac{1}{s^2} + C_2 \frac{1}{s} + C_3 \frac{\ln s}{s} = \mathcal{O} \left( \frac{\ln s}{s} \right)
    }
    $$
    with the constants $C_1, C_2, C_3$ depending on $r, n, d, \| f - \bar{f} \|_F, \| c_a - \bar{c}_a^* \|_F, \Omega$ and $L_{c^*}$.
\end{proof}

\paragraph{Convergence at $1/s^2$ rate}
\begin{proof}[Proof of \cref{prop:piecewise-const-convergence}]
    Let $c_a^*(\omega)$ be a piecewise-constant approximation of $c^*(\omega) = \inf_x f(x, \omega)$ on equidistant grid-points. Then $\|c^* - c^*_a\|_{L^1\Omega} \leq C\frac 1 {s_p}$ where $s_p$ is the number of grid points $\omega_i$.
    Let $$c^*_s(\omega)  = \sum c^*_s(\omega_i) 1_{[\omega_i,\omega_{i+1}]} $$
    where $c^*_s(\omega_i)$ is the best lower bound (resulting from regular SOS) of degree $s$ of $x\mapsto f(x,\omega_i)$.
    Then we have $c_a^*(\omega_i)-c_s^*(\omega_i)$ can be bounded by $$\max_{\omega_i}\|f(\omega_i,\cdot)-\bar f(\omega_i,\cdot)\|_F \left( 1- \left(1-\frac{6r^2}{s^2} \right)^{-n} \right)$$ by \cite{bach_exponential_2023}. Then $$\int_\Omega c^*(\omega)- c_s^*(\omega)\mathrm{d}\omega \leq  \sum_{i} |c_a^*(\omega_i)-c_s^*(\omega_i)||\Delta(\omega_i)| + \|c^* - c_a^* \|_{L^1(\Omega)}.$$
    Using the same bound we get for the first term from the proof of \cref{thm:convergence},
    we can reduce the first term to a $O(1/s^2)$ dependence and we use the theorem on the $L^1$
    convergence of piecewise-constant approximation to 1-periodic trigonometric polynomials from \cite{jackson1930theory} for the second:
    $$\int_\Omega c^*(\omega)- c^*_s(\omega)\mathrm{d}\omega \leq \max_{\omega_i}\|f(\omega_i,\cdot)-\bar f(\omega_i,\cdot)\|_F \left( 1 - \left(1-\frac{6r^2}{s^2}\right)^{-n} \right)|\Omega|+\frac C{s_p}$$
\end{proof}

\subsection{S-SOS for a simple quadratic potential}\label{appendix:simple-potential}

We provide a simple application of S-SOS to a simple quadratic potential that admits a closed-form
solution so as to demonstrate its usage and limitations.

\subsubsection{Analytic solution for the lower bounding function $c^*(\omega)$ with $\omega \sim \text{Uniform}(-1, 1)$}\label{appendix:simple-potential:analytic-solution}
Let $x \in \mathbb{R}$ and $\omega \sim \text{Uniform}(-1, 1)$.
Suppose that we have $$ f(x, \omega) = (x-\omega)^2 + (\omega x)^2 $$

In this case we may explicitly evaluate the exact minimum function
$c^*(\omega) = \inf_x f(x; \omega)$.
Note that
\begin{align*}
f(x; \omega) &= x^2 - 2 \omega x + \omega^2 + \omega^2 x^2
\end{align*}

Explicitly evaluating the zeros of the first derivative we have
$$ \partial_x f(x; \omega) = 2 x^* - 2 \omega + 2 \omega^2 x^* = 0 $$
$$ x^* (1 + \omega^2) = \omega $$
$$ x^* = \frac{\omega}{1 + \omega^2} $$
and, thus, 
$$ c^*(\omega) = \inf_x f(x; \omega) = \frac{\omega^4}{1 + \omega^2}.$$

Note that despite $f(x, \omega)$ being a simple degree-2 SOS polynomial, the tightest lower-bound $c^*(\omega) = \inf_x f(x, \omega)$ is explicitly not polynomial. However, it is algebraic, as it is defined implicitly as the root of the polynomial equation $$c^*(\omega) (1 + \omega^2) - \omega^4 = 0 $$

\subsubsection{Degree-$2s$ S-SOS to find a polynomial lower-bounding function $c^*_{2s}(\omega)$}

Observe that the tightest lower-bounding function $c^*(\omega)$ is not polynomial
even in this simple setting.
However, we can relax the problem to trying to find
$c_{2s} \in \mathcal{P}^{2s}(\Omega)$ to obtain a weaker bound with
$\inf_x f(x, \omega) = c^*(\omega) \geq c_{2s}(\omega)$.

We now proceed with formulating and solving the degree-$2s$ primal S-SOS SDP (\cref{eq:primal-sos-deg}).
We assume that $c_{2s}(\omega)$ is parameterized by a polynomial
of degree $\leq 2s$ in $\omega$.
Observe that this class of functions is not large enough
to contain the true function $c^*(\omega)$.

We choose $s \in \{2, 4\}$ and use the 
standard monomial basis in $x, \omega$,
we have the feature maps $m_2(x, \omega): \mathbb{R}^2 \to \mathbb{R}^6$ and $m_4(x, \omega): \mathbb{R}^2 \to \mathbb{R}^{15}$, since there are $\binom{n+s}{s}$ unique monomials of up to degree-$s$ in $n$ variables.
These assumptions together enable us to explicitly write a SOS SDP in terms of coefficient matching.
Note that we must assume some noise distribution $\nu(\omega)$.
For this section, we present results assuming $\omega \sim \text{Uniform}(-1, 1)$.
We solve the resulting SDP in CVXPY using Legendre quadrature with
$k=5$ zeroes on $[-1, 1]$
to evaluate the objective $\int c(\omega) d\nu(\omega)$.
In fact, $k$ sample points suffice to exactly integrate polynomials of degree
$\leq 2k - 1$.

\begin{figure}
\begin{center}
    \includegraphics[width=0.99\linewidth]{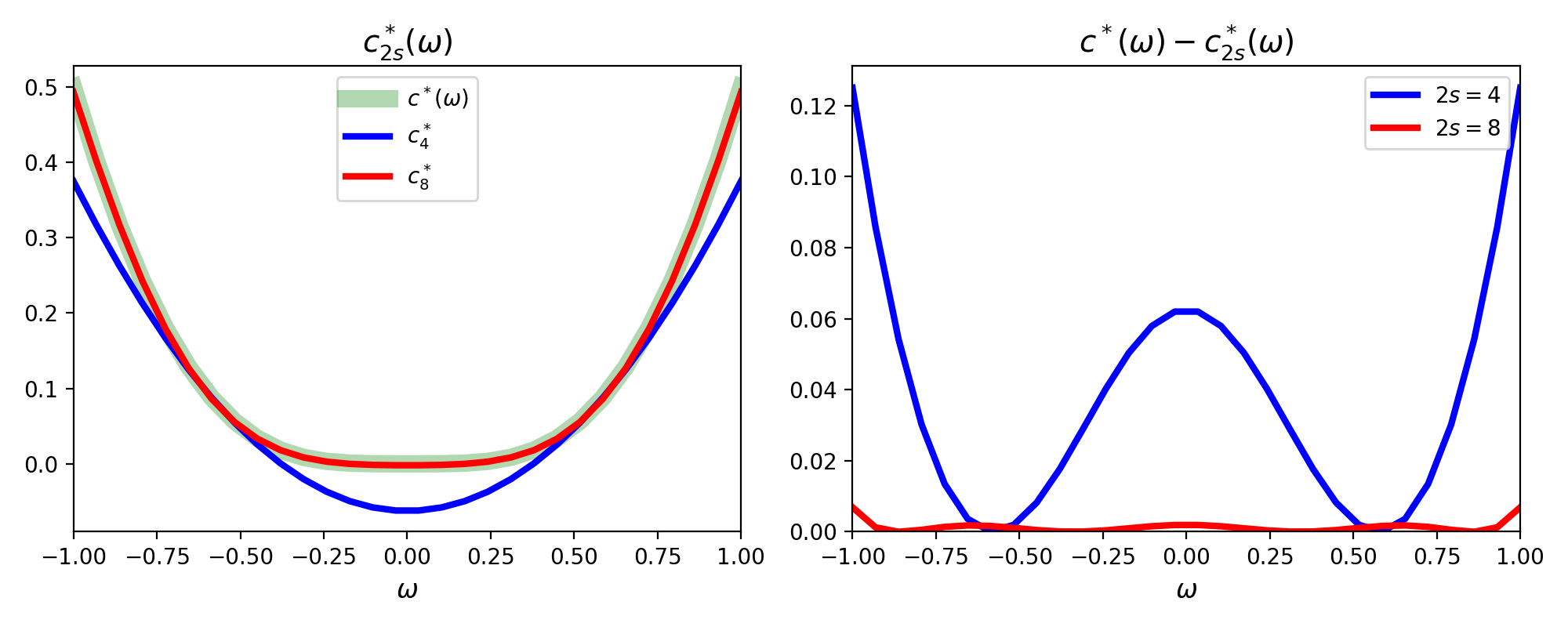}
\end{center}
\caption{Lower bound functions for basis function degree $d=2, 4$ (left)
and the optimality gap to the true lower bound $c^*(\omega) - c^*_{2s}(\omega)$ (right)}
\label{fig:appendix-simple-degree-2-4-gap}
\end{figure}

We solve the SDP for two different levels of the hierarchy, $s=2$ and $s=4$ (producing lower-bound polynomials of degree $4$ and $8$ respectively), and plot the 
lower bound functions $c_{2s}(\omega)$ vs the true lower bound $c^*(\omega) = \omega^4 / (1 + \omega^2)$
as well as the optimality gap to the true lower bound in Fig.\ref{fig:appendix-simple-degree-2-4-gap}.

\subsubsection{Convergence of lower bound as degree $s$ increases}\label{appendix:simple-potential:tightness-as-degree-increases}

To solve the S-SOS SDP in practice, we must 
choose a maximum degree $2s$ for the 
SOS function $m_2(x, \omega)^T W m_2(x, \omega)$
and the lower-bounding function $c(\omega)$, which are both
restricted to be polynomials.
Indeed, a larger $s$ not only increases the dimension
of our basis function $m_s(x, \omega)$ but also the
complexity of the resulting SDP.
We would expect that $d^*_{2s} \to d^*$ as $s \to \infty$,
i.e. the optimal value of the degree-$2s$ S-SOS SDP (\cref{eq:dual-sdp-deg})
converges to that of the ``minimizing distribution'' optimization problem (\cref{eq:dual-sdp}).

In particular, note that in the standard SOS hierarchy
we typically find finite convergence (exact agreement
at some degree $2s^* < \infty$). However, in S-SOS, we thus
far have only a guarantee of asymptotic convergence, as each finite-degree S-SOS SDP solves for a polynomial
approximation to the optimal lower bound
$c^*(\omega) = \inf_{x \in X} f(x, \omega)$. In Figure \ref{fig:obj-value-vs-degree-two-plots}, we illustrate the primal S-SOS SDP objective values
$$ p^*_{2s} = \sup_{c \in \mathcal{P}^{2s}(\Omega)} \int c(\omega) d\nu(\omega) \quad \text{with} \quad f(x, \omega) - c(\omega) \in \mathcal{P}^{2s}_{\text{SOS}}(X \times \Omega) $$
for a given level of the hierarchy (a chosen degree $s$ for the basis $m_s(x, \omega)$) and their convergence
towards the optimal objective value
$$ \int c^*(\omega) d\nu(\omega) = \frac{\pi}{4} - \frac{2}{3} \approx 0.1187 $$
for the simple quadratic potential, assuming
$\nu(\omega) = \frac{1}{2}$ with $\omega \sim \text{Uniform}(-1, 1)$.
We note that in the log-linear plot (right) we have
a ``hinge''-type curve, with a linear decay (in logspace)
and then flattening completely.
This suggests perhaps that in realistic scenarios 
the degree needed to achieve a close approximation
is very low, lower than suggested by our bounds.
The flattening that occurs here is likely due to the numerical tolerance used in our solver (CVXPY/MOSEK),
as increasing the tolerance also increases the asymptotic gap and decreases the degree
at which the gap flattens out.

\subsubsection{Effect of different noise distributions}\label{appendix:simple-potential:different-noise-distributions}

In the previous two sections, we assumed that $\omega \sim \text{Uniform}(-1, 1)$.
This enabled us to solve the primal exactly using Legendre quadrature of polynomials.
Note that in \cref{fig:obj-value-vs-degree-two-plots} we see that the lower-bounding $c^*_2(\omega), c^*_4(\omega)$
for $\omega \sim \text{Uniform}(-1, 1)$ is a smooth polynomial that has curvature (i.e. sign matching
that of the true minimum).
This is actually not guaranteed, as we will see shortly.

\begin{figure}
\begin{center}
    \includegraphics[width=0.80\linewidth]{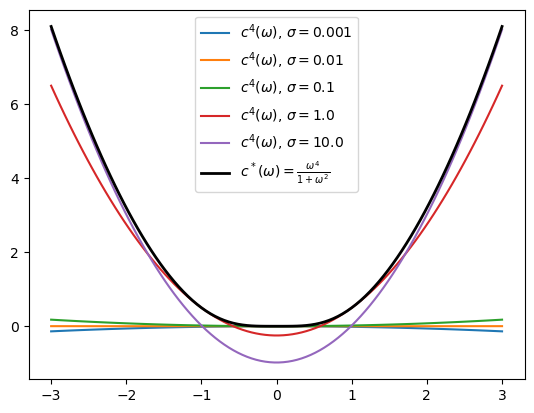}
\end{center}
\caption{
Different lower-bounding functions for degree-4 S-SOS done on the simple quadratic potential
$f(x, \omega) = (x-\omega)^2 + (\omega x)^2$.
The true lower-bounding function $c^*(\omega)$ is plotted in black.
}
\label{fig:simple-univariate-varying-sigma}
\end{figure}

In \cref{fig:simple-univariate-varying-sigma}, we present the lower-bounding functions $c^*_4(\omega)$
achieved by degree-4 S-SOS by solving the dual for $\omega \sim \text{Normal}(0, \sigma^2)$ for varying widths $\sigma$.
We can see that for small $\sigma \ll 1$, the primal solution only cares about the lower-bound accuracy
within a small region of $\omega = 0$, and the lower-bounding curve fails to ``generalize'' effectively
outside the region of consideration.

\subsection{S-SOS for sensor network localization}\label{appendix:snl}

\subsubsection{SDP formulation}\label{appendix:snl:sdp-formulation}

Recall the form of $f(x, \omega)$:
\begin{align*}
f(x, \omega; X, A, r) &= \underbrace{
    \sum_{d_{ij} \in \mathcal{D}_{ss}(r)} ( ||x_i - x_j||_2^2 - d_{ij}(\omega)^2 )^2
}_{\text{sensor-sensor interactions}}
+ \underbrace{
    \sum_{d_{ik} \in \mathcal{D}_{sa}(r)} ( ||x_i - a_k||_2^2 - d_{ik}(\omega)^2 )^2
}_{\text{sensor-anchor interactions}}
\end{align*}

Note that the function $f(x, \omega)$ is exactly a degree-4 SOS polynomial,
so it suffices to choose the degree-2 monomial basis containing $a = \binom{N \ell + d + 2}{2}$ elements as $m_2(x, \omega): \mathbb{R}^{N \ell + d} \to \mathbb{R}^a$.
That is, we have $N$ sensor positions in $\ell$ spatial dimensions and $d$ parameters for a total of
$N \ell + d$ variables.

Let the moment matrix be
$M \in \mathbb{R}^{a \times a}$ with elements defined as
$$ M_{i, j} := \int m_2^{(i)}(x, \omega) m_2^{(j)}(x, \omega) d\mu(x, \omega) $$
for $i, j \in \{ 1, \ldots, a \}$, which fully specifies the minimizing distribution $\mu(x, \omega)$ as in \cref{eq:dual-sdp-deg}.

Our SDP is then of the form
\begin{align*}
d^*_{4} &= \inf_y \sum_\alpha f_\alpha y_\alpha\\ 
              &\text{s.t.} \;\; M(y) \succcurlyeq 0 \\
              &y_{\alpha} = m_{\alpha} \; \forall \; (\alpha, m_{\alpha}) \in \mathcal{M}_{\nu} \\
              &y_{\alpha} = y^*_{\alpha} \; \forall \; (\alpha, y^*_{\alpha}) \in \mathcal{H}
\end{align*}
where $y_{\alpha} = m_{\alpha}$ corresponds to the moment-matching constraints of \cref{eq:dual-sdp-deg}
and $y_{\alpha} = y^*_{\alpha}$ correspond to any possible hard equality constraints required to set the exact position (and uncertainty)
of a sensor $\mathbb{E}[x_i] = x_i^*, \mathbb{E}[x_i^2] - \mathbb{E}[x_i]^2 = 0$ for all $\omega$.
$\mathcal{M}_{\nu}$ represents the $\binom{d+2s}{2s}$ moment-matching constraints necessary for all moments w.r.t. $\omega$
and $\mathcal{H}$ represents the $2 \ell n$ constraints needed to set the exact positions of $n$ known sensor positions in $\mathbb{R}^{\ell}$
(i.e. 1 constraint per sensor and dimension, 2 each for mean and variance).

\subsubsection{Noise types}\label{appendix:snl:noise-types}

In this paper we focus on the linear uniform noise case, as it is a more accurate reflection of measurement noise in true SNL problems. Special robust estimation approaches may be needed to properly handle the outlier noise case.

\begin{itemize}
    \item \textbf{Linear uniform noise}: for a subset of edges we write $d_{ij,k}(\omega) = d_{ij}^* + \epsilon \omega_k$, $\omega_k \sim \text{Uniform}(-1, 1)$,
        and $\epsilon \geq 0$ some noise scale we set.
        The same random variate $\omega_k$ may perturb any number of edges. Otherwise the observed distances are the true distances.
    \item \textbf{Outlier uniform noise}: for a subset of edges we ignore any information in the actual measurement $d_{ij,k} = \omega_k$, $\omega_k \sim \text{Uniform}(0, 2\sqrt{\ell})$ where $\ell$ is the physical dimension of the problem, i.e. $x_i \in \mathbb{R}^{\ell}$.
\end{itemize}

\subsubsection{Algorithms: S-SOS and MCPO}\label{appendix:snl:algorithms}

Here we explicitly formulate MCPO and S-SOS as algorithms.
Let $X = \mathbb{R}^n, \Omega = \mathbb{R}^d$ and use the standard monomial basis.
We write $z = [x_1, \ldots, x_n, \omega_1, \ldots, \omega_d]$.
Our objective is to approximate $c^*(\omega) = \inf_x f(x, \omega)$ for all $\omega$,
with a view towards maximizing $\int c^*(\omega) d\nu(\omega)$ for $\omega$ sampled from
some probability density $\nu(\omega)$.

MCPO (\cref{alg:mcpo}) simply samples $\omega_t$ and finds a set of tuples $(x^*(\omega_t), \omega_t)$
where the optimal minimizer $(x^*(\omega_t, \omega_t)$ is computed using a local optimization
scheme (we use BFGS).

\begin{algorithm}
\caption{Monte Carlo Point Optimization (MCPO)} \label{alg:mcpo}
\begin{algorithmic}[1]
\State \textbf{Input:} Function $f(x; \omega)$, sampler for distribution $\nu(\omega)$, number of samples $T$
\State \textbf{Output:} Approximate integral $\hat{I}$, empirical distribution $p_{\mathcal{D}}(x)$, empirical mean $\mu$, empirical covariance $\Sigma$
\For {$t = 1$ to $T$}
    \State Sample $\omega_t \sim \nu(\omega)$
    \State Find minimizer $x_t = \min_x f(x; \omega_t)$ using BFGS
\EndFor
\State Estimate integral $\hat{I} \approx \frac{1}{T} \sum_{t=1}^{T} f(x_t; \omega_t)$
\State Construct empirical distribution $$ p_{\mathcal{D}}(x) = \frac{1}{T} \sum_{t=1}^{T} \delta(x - x_t) $$
\State Calculate empirical mean $\hat{\mu} = \frac{1}{T} \sum_{t=1}^Tx_t$ and covariance $\hat{\Sigma} = \frac{1}{T-1} \sum_{t=1}^{T} (x_t - \hat{\mu})(x_t - \hat{\mu})^T$.
\end{algorithmic}
\end{algorithm}

S-SOS (\cref{alg:ssos}) via solving the dual (\cref{eq:dual-sdp-deg}) is also detailed below.

\begin{algorithm}
\caption{Stochastic Sum-of-squares (S-SOS), Dual formulation} \label{alg:ssos}
\begin{algorithmic}[1]
\State \textbf{Input:} Maximum basis function degree $s \in \mathbb{Z}_{> 0}$, complete basis function $m_s(x, \omega): \mathbb{R}^{n+d} \to \mathbb{R}^{\binom{n+d+s}{s}}$, function $f(x; \omega): \mathbb{R}^{n+d} \to \mathbb{R}$ represented as a dictionary mapping multi-index $\alpha \in \mathbb{Z}^{n+d}_{\geq 0} \to$ coefficient $f_{\alpha}$, probability density function for $\nu(\omega)$ with known moments $\int \omega^{\alpha} d\nu(\omega) < \infty \; \forall \; ||\alpha||_1 \leq 2s$, any hard equality constraints where we want to set $x_k = x^*_k$ for some $k \in \mathcal{K}$.
\State Let $i1, i2, i4$ be the lexicographically-ordered arrays
    $$\mathbb{Z}_{\geq 0}^{(n+d+1) \times (n+d)}, \mathbb{Z}_{\geq 0}^{\binom{n+d+s}{s} \times (n+d)}, \mathbb{Z}_{\geq 0}^{\binom{n+d+2s}{2s} \times (n+d)}$$
    which correspond to the arrays of multi-indices for all degree-1, degree-$s$, and degree-$2s$ monomials in the variables $z$.
\State Create $M \in \mathbb{R}^{\binom{n+d+s}{s} \times \binom{n+d+s}{s}}$ as a matrix of variables to be estimated.
\State Create $y \in \mathbb{R}^{\binom{n + d + 2s}{2s}}$ as a vector of variables to be estimated, corresponding to the vector of independent moments that fully specifies $M$.
\State Add $M \succcurlyeq 0$ constraint.

\For {$i$ in length($i2$)}
    \For {$j$ in length($i2$)}
        \Comment{Require $M$ to be formed from the elements of $y$.}
        \State Compute $\alpha_{ij} = i2[i] + i2[j]$ as the multi-index corresponding to the sum of the multi-indices $i2[i], i2[j]$.
        \State Add constraint $M_{i, j} = y_{\alpha_{ij}}$.
    \EndFor
\EndFor

\For {each row $\alpha$ in $i4$} \Comment{Require $y_{\alpha}$ moments to equal the known moments of $\omega^{\alpha}$.}
    \If{$\sum_{i=1}^n \alpha_i = 0$}
        \State Add constraint $y_{\alpha} = \int z^{\alpha} d\nu(\omega) = \int \omega^{\alpha[-d:]} d\nu(\omega)$.
    \EndIf
\EndFor

\For {$k$ in $\mathcal{K}$} \Comment{Handle any hard equality constraints in our variables $x$.}
    \State Form multi-index $\alpha_1 \in \mathbb{Z}^{n+d}_{\geq 0}$ where the entry for $x_k$ is set to 1 and everything else is zero.
    \State Form multi-index $\alpha_2 \in \mathbb{Z}^{n+d}_{\geq 0}$ where the entry for $x_k^2$ is set to 1 and everything else is zero.
    \State Add constraint $y_{\alpha_1} = x_k^*$. \Comment{$\mathbb{E}[x_k] = x_k^*$.}
    \State Add constraint $y_{\alpha_2} = (x_k^*)^2$. \Comment{$\text{Var}[x_k] = \mathbb{E}[x_k^2] - \mathbb{E}[x_k]^2 = 0$.}
\EndFor

\State Form the objective to be minimized: $F = \int f(x, \omega) d\mu(x, \omega) = \sum_{\alpha \in i4} f_{\alpha} y_{\alpha}$.
\State Solve SDP where we compute $\inf F$ subject to above constraints.

\State \textbf{Output:} If the problem is feasible (i.e. there exists a degree-$2s$ decomposition of $f$ into $f_{\text{SOS}}$ and $c^*_{2s}(\omega))$, return moment matrix $M \in \mathbb{R}^{\binom{n+d+s}{s} \times \binom{n+d+s}{s}}$, dual objective value $d^*_{2s}$.
Otherwise, terminate and return failed/infeasible SDP solve.
\end{algorithmic}
\end{algorithm}

\subsubsection{Cluster basis hierarchy}\label{appendix:snl:cluster-basis-expanded}

Recall from \cref{section:ssos:variations:cluster-basis} that we defined the cluster basis hierarchy
using body order $b$ and maximum degree per variable $t$.
In this section, we review the additional modifications needed to scale S-SOS for SNL.

In SNL, $f(x, \omega)$ is by design a degree $s=4$ polynomial in $z = [x, \omega]$, with interactions of body order $b=2$ (due to the $(x_i, x_j)$ interactions)
and maximum individual variable degree $t=4$. Written this way, we want to only consider
monomial terms $[x, \omega]^{\alpha}$ with $||\alpha||_1 \leq s$, $||\alpha||_{\infty} \leq 4$, and $||\alpha||_0 \leq 2$.

To sparsify our problem, we start with some $k$-clustering ($k$ clusters, mutually-exclusive) of the sensor set $\mathcal{C} = \{ C_1, \ldots, C_k \}$.
This clustering can be considered as leveraging some kind of ``coarse`` information about which sensors are close to each other.
For example, just looking at the polynomial $f(x, \omega)$ enables us to see which sensors $(i, j)$
must be interacting.

Assume that there is some \emph{a priori} clustering given to us.
We denote $x^{(i)}$ as the subset of the variables restricted to the cluster $C_i$, i.e. $x^{(i)} = \{ x_j: j \in C_i \}$. Moreover, let $G = (V, E)$ be a graph where the vertices $V= \{ 1, \ldots, k \}$ correspond to the $k$ clusters
and the edges $E = \{ (i, j): i, j \in V \}$ correspond to known cluster-cluster interactions.

The SOS part of the 
function $f(x)$ may then be approximated as the sum
of dense intra-cluster interactions and sparse inter-cluster interactions, where the cluster-cluster
interactions are given exactly by edges in the graph $G$:
$$
    m_s(x)^T W m_s(x) \approx
    \sum_{i \in V} m_s(x^{(i)})^T W^{(i)} m_s(x^{(i)})
    + \sum_{(i, j) \in E} m_s(x^{(i)})^T W^{(i, j)} m_s(x^{(j)})
$$
where $W^{(k)}$ are symmetric PSD matrices and $W^{(i, j)}$ are rectangular matrices where we require $W^{(i, j)} = (W^{(j, i)})^T$.
$m_s(x)$ for $x\in \mathbb{R}^n$ here behaves as before and denotes the basis function generated by all $\binom{n+s}{s}$ combinations of monomials with degree $\leq s$.
Notice that this is a strict reduction from the standard Lasserre hierarchy at the same degree $s$, since in general
the standard basis $m_s(x)$ on the full variable set
will contain terms that mix variables from two different
clusters that may not have an edge connecting them.

Efficiency gains in the SDP solve occur when we constrain certain
of the off-diagonal $W^{(i, j)}$ blocks to be zero, i.e. the graph $G$ is sparse in cluster-cluster
interactions.
As we can see from the block decomposition written above, this resembles block sparsity on the matrix $W$.
We may interpret the above scheme as having a hierarchical structure out to depth 2, where we have dense interactions at
the lowest level and sparse interactions aggregating them.
In full generality, the resulting hierarchical sparsity in $W$ may be interpreted as generating a chordal $W$, which
is known to admit certain speed-ups in SDP solvers \cite{vandenberghe_chordal_2017}.

When attempting to solve an SNL problem in the cluster basis
instead of the full basis, we need to throw away terms
in the potential $f(x, \omega)$ that correspond to 
cross-terms that are ``ignored'' by the particular
cluster basis we chose.
The resulting polynomial $\bar{f}(x, \omega)$ has fewer terms and produces a cluster basis SDP that is easier to solve, but generally less accurate due to the sparser
connectivity.

In particular, for the rows in \cref{results:table-ssos-mcpo} that have $N_C > 1$,
we do a $N_C$-means clustering of the ground-truth sensor positions and use those sensor
labels to create our partitioning of the sensors.
We connect every cluster using plus-one $c_i, c_{i+1}$ (including the wrap-around one) connections, so that the cluster-cluster connectivity graph has $N_C$ edges.
We then use this information to throw out observed distances from the set $\mathcal{D}_{ss}$
and from the full basis function $m_2(x, \omega)$.
See our code for complete details.

\subsubsection{Hard equality constraints}\label{appendix:snl:hard-equality-constraints}

The sensor-anchor terms in \cref{eq:snl-soft} are added to make the problem easier,
because by adding them now each sensor no longer needs to rely only on a local neighborhood of sensors
to localize itself, but can also use its position relative to some known anchor.
When we remove them entirely, we need to incorporate hard equality constraints between certain
sensors and known ``anchor'' positions.
This fixes certain known sensors but lets every other sensor be unrooted, defined only relative
to other sensors (and potentially an anchor if it is within the sensing radius).

To deal with the equality constraints
where we set the exact position of a sensor $x_i = x_i^*$,
we solve the dual \cref{eq:dual-sdp-deg} and implement them as equality constraints
on the moment matrix, i.e. for the basis element $m_2(x, \omega)_i = x_i$ we may set $\mathbb{E}[x_i] - x_i^* = M_{0, i} - x_i^* = 0$. Note that we also need to 
set $\text{Var}(x_i) = 0$ so for $m_2(x, \omega)_j = x_i^2$ we add the
equality constraint $\text{Var}(x_i) = \mathbb{E}[x_i^2] - \mathbb{E}[x_i]^2 = M_{0, j} - M_{0, i}^2 = 0$.

\subsubsection{Solution extraction}\label{appendix:snl:solution-extraction}

Once the dual SDP has been solved, we extract the moment matrix $M$ and can easily recover the point and uncertainty estimates
for the sensor positions $\mathbb{E}[x], \text{Var}[x]$ by inspecting the appropriate entries $M_{0, i}$ corresponding to
$m_2(x, \omega)_i = x_i$ and $M_{0, j}$ corresponding to $m_2(x, \omega)_j = x_i^2$.

\subsubsection{Impact of using MCPO with varying numbers of samples $T$}\label{appendix:snl:mcpo-vs-t}

\begin{figure}
\begin{center}
    \includegraphics[width=240px]{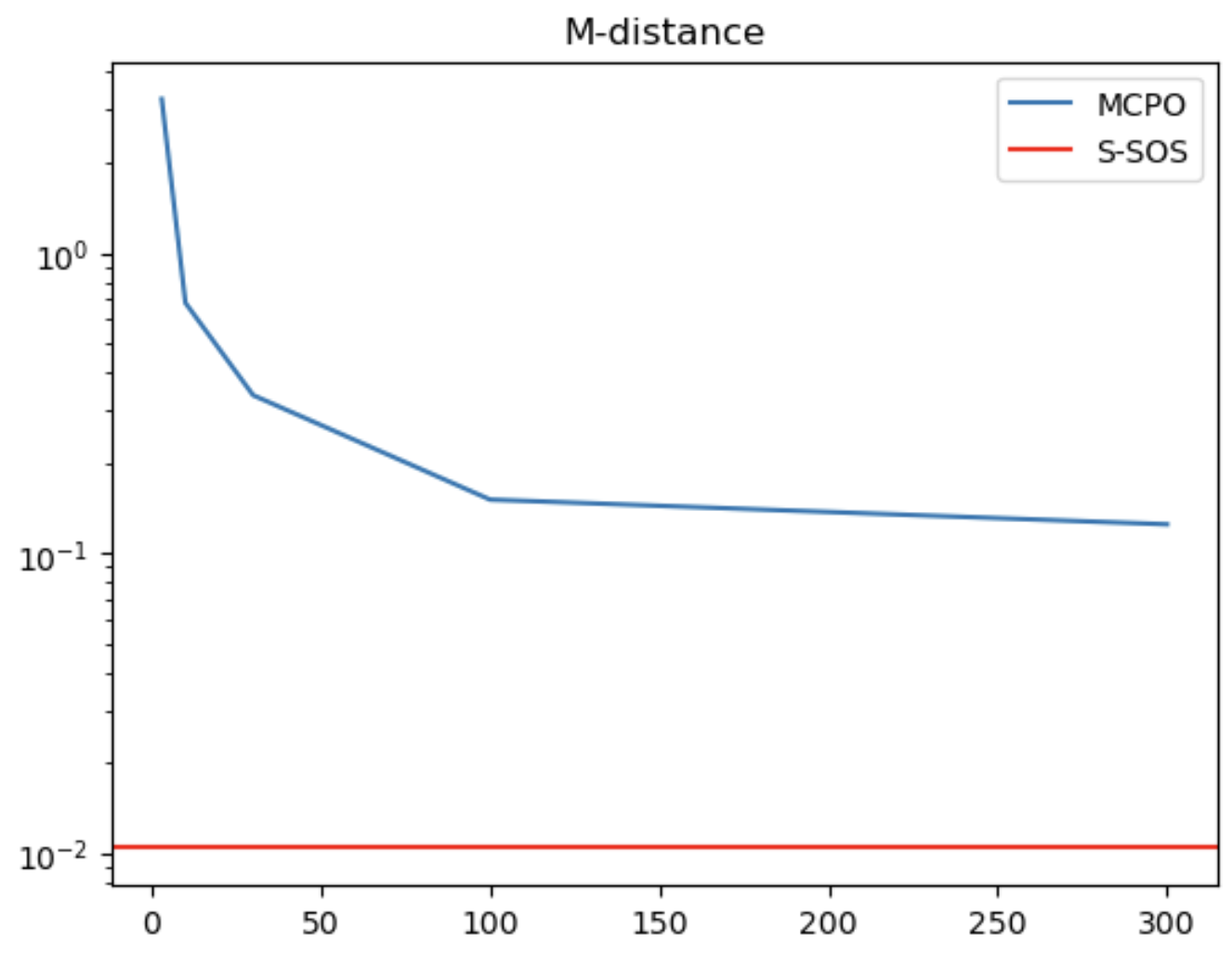}
\end{center}
\caption{Comparison of the performance of MCPO and S-SOS (degree-4) for sensor recovery accuracy in 1D SNL
with varying number of samples $T$ used in the estimate
of empirical $\hat{\mu}, \hat{\Sigma}$.
M-distance is $\delta_M$, our metric for sensor recovery accuracy per \cref{eq:m-distance}.
The problem type here is a $N=5$ sensor, $\ell=1$ spatial dimension, $|\Omega| = 2$ noise variables, $\epsilon = 0.1$ noise scale, $r=3$ sensing radius problem. The full basis is used here for the S-SOS SDP.
}
\label{fig:MCPO-scaling}
\end{figure}

In Figure \ref{fig:MCPO-scaling} we can see how $\delta_M$ varies as we scale the number of samples $T$ used in the MCPO estimate
of the empirical mean/covariance of the recovered solutions. In this particular example, the runtime of the S-SOS
estimate was 0.3 seconds, comparing to 30 seconds for the $T=300$ MCPO point. Despite taking ~100x longer, the MCPO solution recovery still dramatically underperforms
S-SOS in $\delta_M$.
This reflects the poor performance of local optimization methods vs. a global optimization method (when it is available).

\subsubsection{Scalability}\label{appendix:snl:scalability}

The largest 2D SNL experiment we could run had $N=15$ sensors,
$N_C=9$ clusters, and $d=9$ noise parameters.
This generated $N\ell + d = 39$ variables and $820$ basis elements
in the naive $m_2(x, \omega)$ construction, which was reduced to $317$ 
after our application of the cluster basis, giving us
$W, M \in \mathbb{R}^{317 \times 317}$.
A single solve in CVXPY (MOSEK) took 30 minutes on our workstation
(2x Intel Xeon 6130 Gold and 256GB of RAM).
We attempted a run with $N=20$ sensors and $N_C=9$ clusters and $d=9$ noise parameters,
but the process failed due to OOM constraints.
Thus, we report the largest experiment that succeeded.

\end{document}